%
%
%
\input amssym.tex
\input miniltx   
\def\Gin@driver{pdftex.def}
\input graphicx.sty
\resetatcatcode
\newif\ifsect\newif\iffinal
\secttrue\finalfalse
\def\thm #1: #2{\medbreak\noindent{\bf #1:}\if(#2\thmp\else\thmn#2\fi}
\def\thmp #1) { (#1)\thmn{}}
\def\thmn#1#2\par{\enspace{\sl #1#2}\par
        \ifdim\lastskip<\medskipamount \removelastskip\penalty 55\medskip\fi}

\def\pf{\ifdim\lastskip<\smallskipamount \removelastskip\smallskip\fi
        \noindent{\sl Proof\/}:\enspace}
\def\itm#1{\item{\rm #1}\ignorespaces}

\def\bar#1{\overline{#1}}
%
%
\def\Figurascaledeps #1 (#2)(#3){\message{Figura #1}
	\midinsert      
	\centerline{\includegraphics[width=3cm, height=3cm]{#2.pdf}}
	\bigskip
	\centerline{\bf Figure~#1}
	\endinsert}
\newcount\parano
\newcount\eqnumbo
\newcount\thmno
\newcount\versiono
\newcount\remno
\newbox\notaautore
\def\neweqt#1$${\xdef #1{(\number\parano.\number\eqnumbo)}
    \eqno #1$$
    \global \advance \eqnumbo by 1}
\def\newrem #1\par{\global \advance \remno by 1
    \medbreak
{\bf Remark \the\parano.\the\remno:}\enspace #1\par
\ifdim\lastskip<\medskipamount \removelastskip\penalty 55\medskip\fi}
\def\newthmt#1 #2: #3{\xdef #2{\number\parano.\number\thmno}
    \global \advance \thmno by 1
    \medbreak\noindent
    {\bf #1 #2:}\if(#3\thmp\else\thmn#3\fi}
\def\neweqf#1$${\xdef #1{(\number\eqnumbo)}
    \eqno #1$$
    \global \advance \eqnumbo by 1}
\def\newthmf#1 #2: #3{\xdef #2{\number\thmno}
    \global \advance \thmno by 1
    \medbreak\noindent
    {\bf #1 #2:}\if(#3\thmp\else\thmn#3\fi}
\def\forclose#1{\hfil\llap{$#1$}\hfilneg}
\def\newforclose#1{
	\ifsect\xdef #1{(\number\parano.\number\eqnumbo)}\else
	\xdef #1{(\number\eqnumbo)}\fi
	\hfil\llap{$#1$}\hfilneg
	\global \advance \eqnumbo by 1
	\iffinal\else\rsimb#1\fi}
\def\forevery#1#2$${\displaylines{\let\eqno=\forclose
        \let\neweq=\newforclose\hfilneg\rlap{$\qquad\quad\forall#1$}\hfil#2\cr}$$}
\def\noNota #1\par{}
\def\today{\ifcase\month\or
   January\or February\or March\or April\or May\or June\or July\or August\or
   September\or October\or November\or December\fi
   \space\number\year}
\def\inizia{\ifsect\let\neweq=\neweqt\else\let\neweq=\neweqf\fi
\ifsect\let\newthm=\newthmt\else\let\newthm=\newthmf\fi}
\def\bititolo{\empty}
\gdef\begin #1 #2\par{\xdef\titolo{#2}
\ifsect\let\neweq=\neweqt\else\let\neweq=\neweqf\fi
\ifsect\let\newthm=\newthmt\else\let\newthm=\newthmf\fi
\iffinal\let\Nota=\noNota\fi
\centerline{\titlefont\titolo}
\if\bititolo\empty\else\medskip\centerline{\titlefont\bititolo}
\xdef\titolo{\titolo\ \bititolo}\fi
\bigskip
\centerline{\bigfont
\autore \ifvoid\notaautore\else\footnote{${}^1$}{\unhbox\notaautore}\fi}
\bigskip\if\istituto!\centerline{\today}\else
\centerline{\istituto}
\centerline{\indirizzo}
\centerline{\email}
\medskip
\centerline{#1\ \anno}\fi
\bigskip\bigskip
\ifsect\else\global\thmno=1\global\eqnumbo=1\fi}
\def\anno{2011}
\def\raggedleft{\leftskip2cm plus1fill \spaceskip.3333em \xspaceskip.5em
\parindent=0pt\relax}
\font\titlefont=cmssbx10 scaled \magstep1
\font\bigfont=cmr12
\font\eightrm=cmr8
\font\sc=cmcsc10
\font\bbr=msbm10
\font\sbbr=msbm7
\font\ssbbr=msbm5

\font\bfm=cmmib10

\def\ca #1{{\cal #1}}

\nopagenumbers
\binoppenalty=10000
\relpenalty=10000
\newfam\amsfam
\textfont\amsfam=\bbr \scriptfont\amsfam=\sbbr \scriptscriptfont\amsfam=\ssbbr
\newfam\boldifam
\textfont\boldifam=\bfm
\let\de=\partial
\let\eps=\varepsilon
\let\phe=\varphi

\def\Sing{\mathop{\rm Sing}\nolimits}
\def\End{\mathop{\rm End}\nolimits}

\def\Res{\mathop{\rm Res}\nolimits}
\def\Re{\mathop{\rm Re}\nolimits}
\def\Im{\mathop{\rm Im}\nolimits}

\def\id{\mathop{\rm id}\nolimits}

\def\ord{\mathop{\rm ord}\nolimits}

\def\bigoperp{\mathop{\hbox{$\bigcirc\kern-11.8pt\perp$}}\limits}

\mathchardef\void="083F
\mathchardef\ellb="0960
\mathchardef\taub="091C
\def\C{{\mathchoice{\hbox{\bbr C}}{\hbox{\bbr C}}{\hbox{\sbbr C}}
{\hbox{\sbbr C}}}}
\def\R{{\mathchoice{\hbox{\bbr R}}{\hbox{\bbr R}}{\hbox{\sbbr R}}
{\hbox{\sbbr R}}}}
\def\N{{\mathchoice{\hbox{\bbr N}}{\hbox{\bbr N}}{\hbox{\sbbr N}}
{\hbox{\sbbr N}}}}
\def\P{{\mathchoice{\hbox{\bbr P}}{\hbox{\bbr P}}{\hbox{\sbbr P}}
{\hbox{\sbbr P}}}}
\def\Q{{\mathchoice{\hbox{\bbr Q}}{\hbox{\bbr Q}}{\hbox{\sbbr Q}}
{\hbox{\sbbr Q}}}}

\newcount\notitle
\notitle=1
\headline={\ifodd\pageno\rhead\else\lhead\fi}
\def\rhead{\ifnum\pageno=\notitle\iffinal\hfill\else\hfill\tt Version
\the\versiono; \the\day/\the\month/\the\year\fi\else\hfill\eightrm\titolo\hfill
\folio\fi}
\def\lhead{\ifnum\pageno=\notitle\hfill\else\eightrm\folio\hfill\autore\hfill
\fi}
\newbox\bibliobox
\def\setref #1{\setbox\bibliobox=\hbox{[#1]\enspace}
    \parindent=\wd\bibliobox}
\def\biblap#1{\noindent\hang\rlap{[#1]\enspace}\indent\ignorespaces}
\def\art#1 #2: #3! #4! #5 #6 #7-#8 \par{\biblap{#1}#2: {\sl #3\/}.
    #4 {\bf #5} (#6)\if.#7\else, \hbox{#7--#8}\fi.\par\smallskip}
\def\book#1 #2: #3! #4 \par{\biblap{#1}#2: {\bf #3.} #4.\par\smallskip}
\def\coll#1 #2: #3! #4! #5 \par{\biblap{#1}#2: {\sl #3\/}. In {\bf #4,}
#5.\par\smallskip}
\def\pre#1 #2: #3! #4! #5 \par{\biblap{#1}#2: {\sl #3\/}. #4, #5.\par\smallskip}
\def\Bittersweet{\pdfsetcolor{0. 0.75 1. 0.24}} 
\def\Black{\pdfsetcolor{0. 0. 0. 1.}} 
\def\raggedleft{\leftskip2cm plus1fill \spaceskip.3333em \xspaceskip.5em
\parindent=0pt\relax}
\def\Nota #1\par{\medbreak\begingroup\Bittersweet\raggedleft
#1\par\endgroup\Black
\ifdim\lastskip<\medskipamount \removelastskip\penalty 55\medskip\fi}
\newcount\defno\newcount\exno
\def\smallsect #1. #2\par{\bigbreak\noindent{\bf #1.}\enspace{\bf #2}\par
    \global\parano=#1\global\eqnumbo=1\global\thmno=1\global\defno=0\global\remno=0\global\exno=0
    \nobreak\smallskip\nobreak\noindent\message{#2}}
\def\newdef #1\par{\global \advance \defno by 1
    \medbreak
{\bf Definition \the\parano.\the\defno:}\enspace #1\par
\ifdim\lastskip<\medskipamount \removelastskip\penalty 55\medskip\fi}
\def\newex #1\par{\global \advance \exno by 1
    \medbreak
{\bf Example \the\parano.\the\exno:}\enspace #1\par
\ifdim\lastskip<\medskipamount \removelastskip\penalty 55\medskip\fi}
\finaltrue
\versiono=1
%



\secttrue
\def\anno{2015}

\def\autore{Marco Abate\footnote{}{\eightrm 2010 Mathematics Subject Classification: Primary 32H50; 37F75, 37F99.}}
\def\indirizzo{Largo Pontecorvo 5, 56127 Pisa, Italy.\footnote{}{\eightrm Partially supported by FIRB 2012 project ``Geometria Differenziale e Teoria Geometrica delle Funzioni"}}
\def\istituto{Dipartimento di Matematica, Universit\`a di Pisa}
\def\email{E-mail: abate@dm.unipi.it}

\begin {July} Fatou flowers and parabolic curves

{{\narrower{\sc Abstract.} In this survey we collect the main results known up to now (July 2015) regarding possible generalizations to several complex variables of the classical Leau-Fatou flower theorem about holomorphic parabolic dynamics.

}}

\smallsect 1. The original Leau-Fatou flower theorem

In this survey we shall present the known generalizations of the classical Leau-Fatou theorem describing the local holomorphic dynamics about a parabolic point. But let us start with a number of standard definitions.

\newdef A {\sl local $n$-dimensional discrete holomorphic dynamical system} (in short, a {\sl local dynamical system}) is a holomorphic germ $f$ of self-map of a complex $n$-dimensional manifold $M$ at a point $p\in M$ such that
$f(p)=p$; we shall denote by $\End(M,p)$ the set of such germs. 

If $f$, $g$ belongs to $\End(M,p)$ their composition $g\circ f$ is defined as germ in $\End(M,p)$;
in particular, we can consider the sequence $\{f^k\}\subset\End(M,p)$ of {\sl iterates} of $f\in\End(M,p)$, inductively defined by $f^0=\id_M$ and $f^k=f\circ f^{k-1}$ for $k\ge 1$. The aim of local discrete dynamics is exactly the study of the behavior of the sequence of iterates.

\newrem In practice, we shall work with representatives, that is with holomorphic maps $f\colon U\to M$, where $U\subseteq M$ is an open neighborhood of~$p\in U$, such that $f(p)=p$. The fact we are working with germs will be reflected in the freedom we have in taking $U$ as small as needed.
We shall also mostly (but not always) take $M=\C^n$ and $p=O$; indeed a choice of local coordinates $\phe$ for $M$ centered at~$p$ yields an isomorphism $\phe_*\colon\End(M,p)\to\End(\C^n,O)$ preserving the composition by setting $\phe_*(f)=\phe\circ f\circ\phe^{-1}$. 

\newdef Let $f\colon U\to M$ be a representative of a germ in $\End(M,p)$. The {\sl stable set} $K_f\subseteq U$ of~$f$ is the set of points $z\in U$ such that $f^k(z)$ is defined for all $k\in\N$; clearly, $p\in K_f$. If $z\in K_f$, the set $\{f^k(z)\}$ is the {\sl orbit} of~$z$; if $z\in U\setminus K_f$ we shall say that $z$ {\sl escapes.} The stable set depends on the chosen representative, but its germ at~$p$ does not; so we shall freely talk about the stable set of an element of $\End(M,p)$.
An $f$-{\sl invariant} set is a subset $P\subseteq U$ such that $f(P)\subseteq P$; clearly, the stable set is $f$-invariant.

\newdef A local dynamical system $f\in\End(M,p)$ is {\sl parabolic} (and sometimes we shall say that $p$ is a {\sl parabolic fixed point} of~$f$) if $df_p$ is diagonalizable and all its eigenvalues are roots of unity; is {\sl tangent to the identity} if $df_p=\id$. We shall denote by~$\End_1(M,p)$ the set of local dynamical systems tangent to the identity in~$p$.

\newrem If $f\in\End(M,p)$ is parabolic then a suitable iterate $f^q$ is tangent to the identity;
for this reason we shall mostly concentrate on germs tangent to the identity. Furthermore, 
if $f\in\End(M,p)$ is tangent to the identity then $f^{-1}$ is a well-defined germ in $\End(M,p)$ still tangent to the identity.

\newdef The {\sl order} $\ord_p(f)$ of a holomorphic function $f\colon M\to\C$ at~$p\in M$ is the order of vanishing at~$p$, that is the degree of the first non-vanishing term in the Taylor expansion of~$f$ at~$p$ (computed in any set of local coordinates centered at~$p$). The {\sl order}
$\ord_p(F)$ of a holomorphic map $F\colon M\to\C^n$ at~$p\in M$ is the minimum order of its components. 

A germ $f\in\End(\C^n,O)$ can be represented by a $n$-tuple of convergent power series in $n$ variables; collecting terms of the same degree we obtain the homogeneous expansion.

\newdef A {\sl homogeneous map} of degree~$d\ge 1$ is a map $P\colon\C^n\to\C^n$ where $P$ is a $n$-tuple of homogeneous polynomials of degree~$d$ in $n$ variables.
The {\sl homogeneous expansion} of a germ tangent to the identity $f\in\End_1(\C^n,O)$, $f\not\equiv\id_{\C^n}$, is the (unique) series expansion
$$
f(z)=z+P_{\nu+1}(z)+P_{\nu+2}(z)+\cdots
\neweq\equzero
$$
where $P_k$ is a homogeneous map of degree~$k$, and $P_{\nu+1}\not\equiv O$. The number $\nu\ge 1$ is the {\sl order} (or, sometimes, {\sl multiplicity})~$\nu(f)$ of~$f$ at~$O$, and $P_{\nu+1}$ is the {\sl leading term} of~$f$. It is easy to check that the order is invariant under change of coordinates, and thus it can be defined for any germ tangent to the identity $f\in\End_1(M,p)$; we shall denote by $\End_\nu(M,p)$ the set of germs tangent to the identity with order at least~$\nu$. 

In the rest of this section we shall discuss the 1-dimensional case, where the homogeneous expansion reduces to the usual Taylor expansion
$$
f(z)=z+a_{\nu+1}z^{\nu+1}+O(z^{\nu+2})
\neweq\equuno
$$ 
with $a_{\nu+1}\ne 0$. 

\newdef Let $f\in\End_1(\C,0)$ be tangent to the identity given by \equuno. A unit vector $v\in S^1$ is an {\sl attracting} (respectively, {\sl repelling}) {\sl direction} for~$f$ at~$0$ if $a_{\nu+1} v^\nu$ 
is real and negative (respectively, positive). Clearly, there are $\nu$ equally spaced attracting directions, separated by $\nu$ equally spaced repelling directions.

\newex To understand this definition, let us consider the particular case $f(z)=z+az^{\nu+1}$. If $v\in S^1$ is such that $av^\nu>0$ then for every $z\in\R^+v$ we have $f(z)\in\R^+v$ and $|f(z)|>|z|$; in other words, the half-line $\R^+$ is $f$-invariant and repelled from the origin. Conversely, if $v\in S^1$ is such that $av^\nu<0$ then $\R^+v$ is again $f$-invariant
but now $|f(z)|<|z|$ if $z\in\R^+ v$ is small enough; so there is a segment of $\R^+v$ attracted by the origin.

\newrem If $f\in\End_1(\C,0)$ is given by \equuno\ then
$$
f^{-1}(z)=z-a_{\nu+1}z^{\nu+1}+O(z^{\nu+2})\;.
$$
In particular, if $v\in S^1$ is attracting (respectively, repelling) for $f$ then it is repelling (respectively, attracting) for $f^{-1}$, and conversely.

To describe the dynamics of a tangent to the identity germ two more definitions are needed.

\newdef Let~$v\in S^1$ be an attracting direction for
a~$f\in\End_1(\C,0)$ tangent to the identity. The {\sl basin} centered at~$v$ is the set of
points~$z\in K_f\setminus\{0\}$ such that $f^k(z)\to 0$ and $f^k(z)/|f^k(z)|\to
v$ (notice that, up to shrinking the domain of~$f$, we can assume that $f(z)\ne
0$ for all $z\in K_f\setminus\{0\}$). If $z$ belongs to the basin centered
at~$v$, we shall say that the orbit of~$z$ {\sl tends to~$0$ tangent to~$v$.}

A slightly more specialized (but more useful) object is the following: 

\newdef Let $f\in\End_1(\C,0)$ be tangent to the identity. An {\sl
attracting petal} with attracting {\sl central direction}~$v\in S^1$ for
$f$ is an open simply
connected $f$-invariant set $P\subseteq K_f\setminus\{0\}$ with $0\in\de P$ such that a point
$z\in K_f\setminus\{0\}$ belongs to the basin centered at~$v$ if and only if its
orbit intersects~$P$. In other words, the orbit of a point tends to~$0$ tangent
to~$v$ if and only if it is eventually contained in~$P$. A {\sl repelling petal}
(with repelling central direction) is an attracting petal for the inverse
of~$f$.

We can now state the original {\sl Leau-Fatou flower theorem,} describing the dynamics of a one-dimensional tangent to the identity germ in a full neighborhood of the origin (see, e.g., [M] for a modern proof):

\newthm Theorem \flower: (Leau, 1897 [Le]; Fatou, 1919-20 [F1--3]) 
Let
$f\in\End_1(\C,0)$ be tangent to the identity
of order~$\nu\ge 1$. Let $v_1^+,\ldots,v_\nu^+\in
S^1$ be the
$\nu$ attracting directions of~$f$ at the origin, and $v_1^-,\ldots,v_r^-\in
S^1$ the $\nu$ repelling directions. Then:
\itm{(i)} for each attracting (repelling) direction~$v_j^+$ ($v_j^-$)
we can find an attracting (repelling) petal~$P_j^+$ ($P_j^-$) such that the
union of these $2\nu$ petals together with the origin forms a neighborhood
of the origin. Furthermore, the
$2\nu$ petals are arranged cyclically so that two petals intersect if and only if
the angle between their central directions is~$\pi/\nu$.
\itm{(ii)} $K_f\setminus\{0\}$ is the (disjoint) union of the basins
centered at the
$\nu$ attracting directions.
\itm{(iii)} If $B$ is a basin centered at one of the attracting directions, then
there is a function $\chi\colon B\to\C$ such that \indent $\chi\circ
f(z)=\chi(z)+1$ for all $z\in B$. Furthermore, if $P$ is the corresponding petal
constructed in part {\rm (i),} \indent then $\chi|_P$ is a biholomorphism
with an open subset of the complex plane containing a right
half-plane \indent --- and so $f|_P$ is holomorphically
conjugated to the translation $z\mapsto z+1$.

\newdef The function $\chi\colon B\to\C$ constructed in Theorem~\flower.(iii) is a {\sl Fatou coordinate} on the basin~$B$. 

\newrem Up to a linear change of variable, we can assume that $a_{\nu+1}=-1$ in \equuno, so that the attracting directions are the $\nu$-th roots of unity. Given $\delta>0$, 
the set 
$$
D_{\nu,\delta}=\{z\in\C\mid |z^\nu-\delta|<\delta\}
\neweq\equtre
$$
has exactly $\nu$ connected components (each one symmetric with respect to a different $\nu$-th root of unity), and it turns out that when $\delta>0$ is small enough these components can be taken as attracting petals for~$f$ --- even though to cover a neighborhood of the origin one needs slightly larger petals. The components of $D_{\nu,\delta}$ are distributed as petals in a flower; this is the reason why Theorem~\flower\ is called ``flower theorem".


So the union of attracting and repelling petals gives a pointed neighborhood of the origin, and 
the dynamics of $f$ on each petal is conjugated to a translation via a Fatou coordinate. The relationships between different Fatou coordinates is the key to \'Ecalle-Voronin holomorphic classification of parabolic germs (see, e.g., [A4] and references therein for a concise introduction to \'Ecalle-Voronin invariants), which is however outside of the scope of this survey. We end this section with the statement of the Leau-Fatou flower theorem for general parabolic germs: 

\newthm Theorem \pflower: (Leau, 1897 [Le]; Fatou, 1919-20 [F1--3]) 
Let $f\in\End(\C,0)$ be of the form
$f(z)=\lambda z+O(z^2)$,
where $\lambda\in S^1$ is a
primitive root of the unity of order $q$. Assume that
$f^q\not\equiv\id$. Then there exists $\mu\ge 1$ such that $f^q$ has order~$q\mu$, and $f$ acts on the
attracting (respectively, repelling) petals of $f^q$ as a permutation
composed by $\mu$ disjoint cycles. Finally, $K_f=K_{f^q}$.

In the subsequent sections we shall discuss known generalizations of Theorem~\flower\ to several variables.

\smallsect 2. \'Ecalle-Hakim theory

From now on we shall work in dimension $n\ge 2$. So let $f\in\End_1(\C^n,O)$ be tangent to the identity; we would like to find a multidimensional version of the petals of Theorem~\flower.

If $f$ had a non-trivial one-dimensional $f$-invariant curve passing through the origin, that is an injective holomorphic map $\psi\colon\Delta\to\C^n$, where $\Delta\subset\C$ is a neighborhood of the origin, such that $\psi(0)=O$,
$\psi'(0)\ne O$ and $f\bigl(\psi(\Delta)\bigr)\subseteq\psi(\Delta)$ with  $f|_{\psi(\Delta)}\not\equiv\id$, we could apply Leau-Fatou flower theorem to~$f|_{\psi(\Delta)}$ obtaining a one-dimensional Fatou flower for~$f$ inside
the invariant curve. In particular, if $z^o\in\psi(\Delta)$ belongs to an attractive petal, we 
would have $f^k(z^o)\to O$ and $[f^k(z^o)]\to[\psi'(0)]$, where $[\cdot]\colon\C^n\setminus\{O\}\to\P^{n-1}(\C)$ is the canonical projection. The first observation we can make is that then $[\psi'(0)]$ cannot be any direction in $\P^{n-1}(\C)$. Indeed:

\newthm Proposition \Hakimz: ([H2]) Let $f(z)=z+P_{\nu+1}(z)+\cdots\in\End_1(\C^n,O)$ be tangent to the identity of order $\nu\ge 1$. Assume there is $z^o\in K_f$ such that $f^k(z^o)\to O$
and $[f^k(z^o)]\to[v]\in\P^{n-1}(\C)$. Then $P_{\nu+1}(v)=\lambda v$ for some $\lambda\in\C$.

\newdef Let $P\colon\C^n\to\C^n$ be a homogeneous map.
A direction $[v]\in\P^{n-1}(\C)$ is {\sl characteristic} for~$P$ if $P(v)=\lambda v$ for some $\lambda\in\C$. Furthermore, we shall say that $[v]$ is {\sl degenerate} if~$P(v)=O$, and {\sl non-degenerate} otherwise. 

\newrem From now on, given $f\in\End_1(\C^n,O)$ tangent to the identity of order~$\nu\ge 1$, every notion/object/concept introduced for its leading term $P_{\nu+1}$ will be introduced also for~$f$; for instance, 
a (degenerate/non-degenerate) characteristic direction  for~$P_{\nu+1}$ will also be a (degenerate/non-degenerate) characteristic direction for~$f$.

\newrem If $f\in\End_1(\C^n,O)$ is given by \equzero, then $f^{-1}\in\End_1(\C^n,O)$ is given by 
$$
f^{-1}(z)=z-P_{\nu+1}(z)+\cdots\;.
$$
In particular, $f$ and $f^{-1}$ have the same (degenerate/non-degenerate) characteristic directions.

\newrem If $\psi\colon\Delta\to\C^n$ is a one-dimensional curve with $\psi(0)=O$ and $\psi'(0)\ne O$ such that $f|_{\psi(\Delta)}\equiv\id$, it is easy to see that $[\psi'(0)]$ must be a degenerate characteristic direction for~$f$.

So if we have an $f$-invariant one-dimensional curve~$\psi$ through the origin then $[\psi'(0)]$ must be a characteristic direction. However, in general the converse is false: there are non-degenerate characteristic directions which are not tangent to any $f$-invariant curve passing through the origin.

\newex ([H2]) Let $f\in\End(\C^2,O)$ be given by 
$$
f(z,w)=\left({z\over 1+z}, w
+z^2\right)\;,
$$
so that $f$ is tangent to the identity of order~1, and 
$P_2(z,w)=(-z^2, 
z^2)$. In particular, $f$ has a degenerate characteristic direction $[0:1]$ and a non-degenerate characteristic direction $[v]=[1:
-1]$. The degenerate characteristic direction is tangent to the curve $\{z=0\}$, which is pointwise fixed by~$f$, in accord with Remark~2.3. We claim that no $f$-invariant curve can be tangent to~$[v]$.\hfill\break\indent
Assume, by contradiction, that we have an $f$-invariant curve $\psi\colon\Delta\to\C^2$ with
$\psi(0)=O$ and $[\psi'(0)]=[v]$. Without loss of generality, we can assume that $\psi(\zeta)=\bigl(\zeta,u(\zeta)\bigr)$ with $u\in\End(\C,0)$. Then the condition of $f$-invariance becomes
$f_2\bigl(\zeta,u(\zeta)\bigr)=u\bigl(f_1(\zeta,u(\zeta)\bigr)\bigr)$, that is
$$
u(\zeta)
+\zeta^2=u\left({\zeta\over 1+\zeta}\right)\;.
\neweq\eqduno
$$
Put $g(\zeta)=\zeta/(1+\zeta)$, so that $g^k(\zeta)=\zeta/(1+k\zeta)$; in particular, $g^k(\zeta)\to 0$
for all $\zeta\in\C\setminus\{-{1\over n}\mid n\in\N^*\}$. This means that by using \eqduno\ we can extend $u$ to $\C\setminus\{-{1\over n}\mid n\in\N^*\}$ by setting
$$
u(\zeta)=u\bigl(g^k(\zeta)\bigr)-\sum_{j=0}^{k-1}[g^j(\zeta)]^2
$$
where $k\in\N$ is chosen so that $g^k(\zeta)\in\Delta$. Analogously, 
\eqduno\ implies that for $|\zeta|$ small enough one has
$$
u\bigl(g^{-1}(\zeta)\bigr)
+\bigl(g^{-1}(\zeta)\bigr)^2=u(\zeta)\;;
$$
so we can use this relation to extend $u$ to all of~$\C$, and then to $\P^1(\C)$, because $g^{-1}(\infty)=-1$. So $u$ is a holomorphic function defined on~$\P^1(\C)$, that is a constant; but no constant can satisfy \eqduno, contradiction.

\newrem Rib\'on [R] has given examples of germs having no holomorphic invariant curves at all. For instance, this is the case for
germs of the form $f(z,w)=(z+w^2, w+z^2+\lambda z^5)$ for all $\lambda\in\C$ outside a polar Borel set.

The first important theorem we would like to quote is due to \'Ecalle [E] and Hakim [H2], and it says that we do always have a Fatou flower tangent to a non-degenerate characteristic direction, even when there are no invariant complex curves containing the origin in their relative interior. To state it, we need to define what is the correct multidimensional notion of petal. 

\newdef A {\sl parabolic curve} for $f\in\End_1(\C^n,O)$
tangent to the identity is an injective holomorphic map
$\phe\colon D\to\C^n\setminus\{O\}$ satisfying the following properties:
\smallskip
\item{(a)} \ $D$ is a simply connected domain in~$\C$ with $0\in\de D$;
\item{(b)} \ $\phe$ is continuous at the origin, and $\phe(0)=O$;
\item{(c)} \ $\phe(D)$ is $f$-invariant, and $(f|_{\phe(D)})^k\to O$
uniformly on compact subsets as $k\to+\infty$.
\smallskip
\noindent Furthermore, if $[\phe(\zeta)]\to[v]$ in $\P^{n-1}(\C)$ as~$\zeta\to
0$ in~$D$, we shall say that the parabolic curve $\phe$ is {\sl tangent} to
the direction~$[v]\in\P^{n-1}(\C)$. Finally, a {\sl Fatou flower} with $\nu$ {\sl petals} tangent to 
a direction~$[v]$ is a holomorphic map $\Phi\colon D_{\nu,\delta}\to\C$, where $D_{\nu,\delta}$ is given by \equtre, such that $\Phi$ restricted to any connected component of~$D_{\nu,\delta}$ is a parabolic curve tangent to~$[v]$, a {\sl petal} of the Fatou flower. If $\nu$ is the order of~$f$ then we shall talk of a Fatou flower for~$f$ without mentioning the number of petals.

Then \'Ecalle, using his resurgence theory (see, e.g., [S] for an introduction to \'Ecalle's resurgence theory in one dimension), and Hakim, using more classical methods, have proved the following result (see also [W]):

\newthm Theorem \EcalleHakim: (\'Ecalle, 1985 [E]; Hakim, 1997 [H2, 3]) 
Let $f\in\End(\C^n,O)$ be tangent to the identity,
and $[v]\in\P^{n-1}(\C)$ a non-degenerate characteristic
direction for~$f$. Then there exists (at least) one Fatou flower tangent to~$[v]$. Furthermore, for every petal $\phe\colon\Delta\to\C^n$ of the Fatou flower there exists a injective holomorphic map $\chi\colon\phe(\Delta)\to\C$ such that $\chi\bigl(f(z)\bigr)=\chi(z)+1$ for all $z\in\phe(\Delta)$.

\newdef The function $\chi$ constructed in the previous theorem is a {\sl Hakim-Fatou coordinate}.

\newrem A characteristic direction is a {\it complex} direction, not a real one; so it should not be confused with the attracting/repelling directions of Theorem~\flower. All petals of a Fatou flower are tangent to the same characteristic direction, but each petal is tangent to a different real direction inside the same complex (characteristic) direction. In particular, Fatou flowers of~$f$ and~$f^{-1}$ are tangent to the same characteristic directions (see Remark~2.2) but the corresponding petals are tangent to different real directions, as in Theorem~\flower.


In particular there exist parabolic curves tangent to $[1:-1]$ for the system of Example~2.1
even though there are no invariant curves passing through the origin tangent to that direction.

Parabolic curves are one-dimensional objects in an $n$-dimensional space; it is natural to wonder about the existence of higher dimensional invariant subsets. A sufficient condition for their existence has been given by Hakim; to state it we need to introduce another definition.

\newdef Let $[v]\in\P^{n-1}(\C)$ be a non-degenerate characteristic direction for a homogeneous map $P\colon\C^n\to\C^n$ of degree $\nu+1\ge 2$; in particular, $[v]$ is a fixed point for the meromorphic self-map $[P]$ of~$\P^{n-1}(\C)$ induced by~$P$. The {\sl directors} of~$P$ in~$[v]$ are the eigenvalues $\alpha_1,\ldots,\alpha_{n-1}\in\C$ of the linear operator
 $$
{1\over\nu}\left(d[P]_{[v]}-\id\right)\colon T_{[v]}\P^{n-1}(\C)\to
  T_{[v]}\P^{n-1}(\C)\;.
 $$ 
As usual, if $f\in\End_1(\C^n,O)$ is of the form \equuno, then the {\sl directors} of~$f$ in a non-degenerate characteristic direction $[v]$ are the directors of~$P_{\nu+1}$ in~$[v]$.

\newrem Definition~2.4 is equivalent to the original definition used by Hakim (see, e.g., [ArR]). Furthermore, in dimension~2 if $[v]=[1:0]$ is a non-degenerate characteristic direction of $P=(P_1,P_2)$ we have
$P_1(1,0)\ne 0$, $P_2(1,0)=0$ and the director is given by
$$
{1\over\nu}\left.{d\over d\zeta}{P_2(1,\zeta)-\zeta P_1(1,\zeta)\over P_1(1,\zeta)}\right|_{\zeta=0}
={1\over\nu}\left[{{\de P_2\over\de z_2}(1,0)\over P_1(1,0)}-1\right]\;.
$$ 

\newrem Recalling Remark~2.2 one sees that a germ $f\in\End_1(\C^n,O)$ tangent to the identity and its inverse~$f^{-1}$ have the same directors at their non-degenerate characteristic directions. 

\newrem The proof of Theorem~\EcalleHakim\ becomes simpler when no director is of the form~$k\over\nu$ with $k\in\N^*$; furthermore, in this case the parabolic curves 
enjoy additional properties (in the terminology of [AT1] they are {\sl robust;} see also [Ro3]). 

\newdef A {\sl parabolic manifold} for a germ $f\in\End_1(\C^n,O)$ tangent to the identity is an $f$-invariant complex submanifold $M\subset\C^n\setminus\{O\}$ with $O\in\de M$ such that $f^k(z)\to O$ for all $z\in M$. A {\sl parabolic domain} is a parabolic manifold of dimension~$n$. We shall say that $M$ is {\sl attached} to the characteristic direction $[v]\in\P^{n-1}(\C)$ if furthermore
$[f^k(z)]\to[v]$ for all $z\in M$.

Then Hakim has proved (see also [ArR] for the details of the proof) the following theorem:

\newthm Theorem \Hakimb: (Hakim, 1997 [H3])
Let $f\in\End_1(\C^n,O)$ be tangent to the identity of order~$\nu\ge 1$.
Let
$[v]\in\P^{n-1}(\C)$ be a non-degenerate characteristic direction, with directors
$\alpha_1,\ldots,\alpha_{n-1}\in\C$.
Furthermore, assume that $\Re\alpha_1,\ldots,\Re\alpha_d>0$ and
$\Re\alpha_{d+1},\ldots,\Re\alpha_{n-1}\le 0$ for a suitable $d\ge 0$. Then:
\smallskip
\itm{(i)} There exist (at least) $\nu$ parabolic $(d+1)$-manifolds~$M_1,
\dots, M_\nu$ of~$\C^n$ attached to~$[v]$;
\itm{(ii)} $f|_{M_j}$ is holomorphically conjugated to the translation
$\tau(w_0,w_1,\ldots,w_d)=(w_0+1,w_1,\ldots,w_d)$ defined \indent on a
suitable right half-space in~$\C^{d+1}$. 

\newrem
In particular, if all the directors of~$[v]$
have positive real part, there is at least one parabolic domain.
However, the condition given by Theorem~\Hakimb\ is not necessary for the
existence of parabolic domains; see [Ri1], [Us], [AT3] for examples, and
[Ro8] for conditions ensuring the existence of a parabolic domain when some directors
have positive real part and all the others are equal to zero. Moreover, Lapan [L1] has proved that if $n=2$ and $f$ has a unique characteristic direction~$[v]$ which is non degenerate then there exists a parabolic domain attached to~$[v]$ even though the director is necessarily~0. 

Two natural questions now are: how many characteristic directions are there? Does there always exist a non-degenerate characteristic direction? To answer the first question, we need to introduce the notion of multiplicity of a characteristic direction. To do so, notice that $[v]=[v_1:\cdots:v_n]\in\P^{n-1}(\C)$ is a characteristic direction for the homogeneous map $P=(P_1,\ldots,P_n)$ if and only if $v_hP_k(v)-v_kP_h(v)=0$ for all $h$,~$k=1,\ldots,n$. In particular, the set of characteristic directions of~$P$ is an algebraic subvariety of~$\P^{n-1}(\C)$.

\newdef If the maximal dimension of the irreducible components of the subvariety of characteristic directions of a homogeneous map $P\colon\C^n\to\C^n$ is~$k$, we shall say that $P$ is {\sl $k$-dicritical}; if $k=n$ we shall say that $P$ is {\sl dicritical;} if $k=0$ we shall say that $P$ is {\sl non-dicritical}. 

\newrem A homogeneous map $P\colon\C^n\to\C^n$ of degree~$d$ is dicritical if and only if $P(z)=p(z)z$ for some homogenous polynomial $p\colon\C^n\to C$ of degree~$d-1$. In particular, the degenerate characteristic directions are the zeroes of the polynomial~$p$.

In the non-dicritical case we can count the number of characteristic directions, using a suitable multiplicity. 

\newdef Let $[v]=[v_1:\cdots:v_n]\in\P^{n-1}(\C)$ be a characteristic direction of a homogeneous map $P=(P_1,\ldots,P_n)\colon\C^n\to\C^n$. Choose $1\le j_0\le n$ so that $v_{j_0}\ne 0$. The {\sl multiplicity}~$\mu_P([v])$ of~$[v]$ is the local intersection multiplicity at~$[v]$ in $\P^{n-1}(\C)$ of the polynomials $z_{j_0}P_j-z_jP_{j_0}$ with $j\ne j_0$ if $[v]$ is an isolated characteristic direction; it is $+\infty$ if~$[v]$ is not isolated. 

\newrem The local intersection multiplicity $I(p_1,\ldots,p_k;z^o)$ of a set $\{p_1,\ldots, p_k\}$ of holomorphic functions at a point~$z^o\in\C^n$ can be defined (see, e.g., [GH]) as 
$$
I(p_1,\ldots,p_k;z^o)=\dim \ca O_{n,z^o}/(p_1,\ldots,p_k)\;,
$$
where $\ca O_{n,z^o}$ is the local ring of germs of holomorphic functions at~$z^o$, and the dimension is as vector space. It is easy to check that the definition of multiplicity of a characteristic direction does not depend on the index~$j_0$ chosen. Furthermore, since the local intersection multiplicity is invariant under change of coordinates, we can use local charts to compute the local intersection multiplicity on complex manifolds.

\newrem When $n=2$, the multiplicity of $[v]=[1:v_2]$ as characteristic direction of $P=(P_1,P_2)$
is the order of vanishing at $t=v_2$ of $P_2(1,t)-tP_1(1,t)$; analogously, the multiplicity of $[0:1]$ is the order of vanishing at $t=0$ of $P_1(t,1)-tP_2(t,1)$. 

Then we have the following result (see, e.g., [AT1]):

\newthm Proposition \ncd: Let $P\colon\C^n\to\C^n$ be a non-dicritical homogeneous map of degree~$\nu+1\ge 2$. Then $P$ has exactly 
$$
{1\over\nu}\bigl((\nu+1)^n-1\bigr)=\sum_{j=0}^{n-1} {n\choose j+1}\nu^j
$$
characteristic directions, counted according to their multiplicity.

In particular, when $n=2$ then a homogeneous map of degree~$\nu+1$ either is dicritical (and all directions are characteristic) or has exactly $\nu+2$ characteristic directions. But all of them can be degenerate; an example is the following (but it is easy to build infinitely many others). 

\newex Let $P(z,w)=(z^2w+zw^2,zw^2)$. Then the characteristic directions of~$P$ are $[1:0]$ and $[0:1]$, both degenerate. Using Remark~2.12, we see that $\mu_P([1:0])=3$ and $\mu_P([0:1])=1$. 

So we cannot apply Theorem~\EcalleHakim\ to any germ of the form $f(z)=z+P(z)+\cdots$ when
$P$ is given by Example~2.2. However, as soon as the higher order terms are chosen so that the origin is an isolated fixed point then $f$ does have parabolic curves:

\newthm Theorem \Abate: (Abate, 2001 [A2]) Let $f\in\End_1(\C^2,O)$ be tangent to the identity such that $O$ is an isolated fixed point. Then $f$ admits at least one Fatou flower tangent to some characteristic direction.

In the next section we shall explain why this theorem holds, we shall give more general statements, and we shall give an example (Example 3.1) showing the necessity of the hypothesis that the origin is an isolated fixed point.

\smallsect 3. Blow-ups, indices and Fatou flowers

In the previous section we saw that for studying the dynamics of a germ tangent to the identity it is useful to consider the tangent directions at the fixed point. A useful way for dealing with tangent directions consists, roughly speaking, in replacing the fixed point by the projective space of the tangent directions, in such a way that the new space is still a complex manifolds, where the tangent directions at the original fixed point are now points. We refer to, e.g., [GH] or [A1] for a precise  description of this construction; here we shall limit ourselves to explain how to work with it.

\newdef Let $M$ be a complex $n$-dimensional manifold, and $p\in M$. The {\sl blow-up} of~$M$ of {\sl center}~$p$ is a complex $n$-dimensional manifold $\tilde M$ equipped with a surjective holomorphic map $\pi\colon \tilde M\to M$ such that 
\smallskip
\item{(i)} $E=\pi^{-1}(p)$ is a compact submanifold of~$\tilde M$, the {\sl exceptional divisor} of the blow-up, biholomorphic to $\P(T_pM)$;
\item{(ii)} $\pi|_{\tilde M\setminus E}\colon \tilde M\setminus E\to M\setminus\{p\}$ is a biholomorphism.

Let us describe the construction for $(M,p)=(\C^n,O)$; using local charts one can repeat the construction for any manifold. As a set, $\widetilde{\C^n}$ is the disjoint union of $\C^n\setminus\{O\}$ and $E=\P^{n-1}(\C)$; we shall define a manifold structure using charts.
For $j=1,\ldots,n$ let $U'_j=\{[v_1:\cdots:v_n]\in\P^{n-1}(\C)\mid v_j\ne 0\}$, $U''_j=\{w\in\C^n\mid w_j\ne 0\}$ and $\tilde U_j=U'_j\cup U''_j\subset\widetilde{\C^n}$. Define $\chi_j\colon\tilde U_j\to\C^n$ by setting
$$
\chi_j(q)=\cases{ \left({v_1\over v_j},\ldots,{v_{j-1}\over v_j},0,{v_{j+1}\over v_j},\ldots,{v_n\over v_j}\right)&if $q=[v_1:\cdots:v_n]\in U'_j$,\cr 
\left({w_1\over w_j},\ldots, {w_{j-1}\over w_j},w_j,{w_{j+1}\over w_j},\ldots, {w_n\over w_j}\right)& if $q=(w_1,\ldots,w_n)\in U''_j$.\cr}
$$
We have
$$
\chi^{-1}_j(w)=\cases{[w_1:\cdots:w_{j-1}:1:w_{j+1}:\cdots:w_n]&if $w_j=0$,\cr
(w_jw_,\ldots,w_jw_{j-1},w_j,w_jw_{j+1},\ldots,w_jw_n)& if $w_j\ne 0$,\cr}
$$
and it is easy to check that $\{(\tilde U_1,\chi_1),\ldots,(\tilde U_n,\chi_n)\}$ is an atlas for $\widetilde{\C^n}$, with $\chi_j([0:\cdots:1:\cdots:0])=O$ and $\chi_j\bigl(\tilde U_j\cap\P^{n-1}(\C)\bigr)=\{w_j=0\}\subset\C^n$. We can then define the projection $\pi\colon\widetilde{\C^n}\to\C^n$
in coordinates by setting
$$
\pi\circ\chi_j^{-1}(w)=(w_1w_j,\ldots,w_{j-1}w_j,w_j,w_{j+1}w_j,\ldots,w_nw_j)\;;
$$
it is easy to check that $\pi$ is well-defined, that $\pi^{-1}(O)=\P^{n-1}(\C)$ and that 
$\pi$ induces a biholomorphism between $\widetilde{\C^n}\setminus \P^{n-1}(\C)$ and $\C^n\setminus\{O\}$. 
Notice furthermore that $\widetilde{\C^n}$ has a canonical structure of line bundle over~$\P^{n-1}(\C)$ given by the projection $\tilde\pi\colon\widetilde{\C^n}\to\P^{n-1}(\C)$ defined by
$$
\tilde\pi(q)=\cases{[v]&if $q=[v]\in\P^{n-1}(\C)$,\cr
[w_1:\cdots:w_n]& if $q=w\in\C^n\setminus\{O\}$;\cr}
$$
the fiber over $[v]\in\P^{n-1}(\C)$ is given by the line $\C v\subset\C^n$. 

Two more definitions we shall need later on: 

\newdef Let $\pi\colon\tilde M\to M$ be the blow-up of a complex manifold~$M$ at~$p\in M$. Given a subset $S\subset M$, the {\sl full} (or {\sl total}) {\sl transform} of~$S$ is~$\pi^{-1}(S)$, whereas the {\sl strict transform} of~$S$ is the closure in~$\tilde M$ of~$\pi^{-1}(S\setminus\{O\})$. 

Clearly, the full and the strict transform coincide if $p\notin S$; if $p\in S$ then the full transform is the union of the strict transform and the exceptional divisor. Furthermore,
if $S$ is a submanifold at~$p$ then its strict transform is $(S\setminus\{p\})\cup\P(T_pS)$. 

\newdef Let $f\in\End(M,p)$ be a germ such that $df_p$ is invertible. Choose a representative $(U,f)$ of the germ such that $f$ is injective in~$U$. Then the {\sl blow-up} of~$f$ is 
the map $\tilde f\colon\pi^{-1}(U)\to\tilde M$ defined by
$$
\tilde f(q)=\cases{[df_p(v)]& if $q=[v]\in E=\P(T_pM)$,\cr 
f(w)& if $q=w\in U\setminus\{p\}$.\cr}
$$
In this way we get a germ about the exceptional divisor of a holomorphic self-map of the blow-up, given by the differential of~$f$ along the exceptional divisor and by $f$ itself elsewhere, satisfying 
$\pi\circ\tilde f=f\circ\pi$. In particular, $K_{\tilde f}=\pi^{-1}(K_f)=\bigl(K_f\setminus\{O\}\bigr)\cup
E$, and to study the dynamics of $\tilde f$ in a neighborhood of the exceptional divisor is equivalent to studying the dynamics of~$f$ in a neighborhood of~$p$. 

If $f\in\End_1(\C^n,O)$ is tangent to the identity, its blow-up $\tilde f$ in the chart $(U_1,\chi_1)$ is given by
$$
\chi_1\circ\tilde f\circ\chi_1^{-1}(w)=\cases{w& if $w_1=0$,\cr
\left(f_1(w_1,w_1w_2,\ldots,w_1w_n),{f_2(w_1,w_1w_2,\ldots,w_1w_n)\over f_1(w_1,w_1w_2,\ldots,w_1w_n)},\ldots,{f_n(w_1,w_1w_2,\ldots,w_1w_n)\over f_1(w_1,w_1w_2,\ldots,w_1w_n)}\right)& if $w_1\ne 0$;\cr}
$$
similar formulas hold in the other charts. In particular, writing $w=(w_1,w')$ and $\chi_1\circ\tilde f\circ\chi_1^{-1}=(\tilde f_1,\ldots,\tilde f_n)$, if $f$ is tangent to the identity of order $\nu\ge 1$ and leading term $P_{\nu+1}=(P_{\nu+1,1},\ldots, P_{\nu+1,n})$, we get
$$
\cases{\tilde f_1(w)=w_1+w_1^{\nu+1}P_{\nu+1,1}(1,w')+O(w_1^{\nu+2})\;,&\cr
\tilde f_j(w)=w_j+w_1^\nu\bigl(P_{\nu+1,j}(1.w')-w_jP_{\nu+1,1}(1,w')\bigr)+O(w_1^{\nu+1})&
if $j\ne 1$.\cr}
\neweq\eqtuno
$$ 
It follows immediately that:
\smallskip
\item{--} if $\nu\ge 2$ then $\tilde f$ is tangent to the identity in all points of the exceptional divisor;
\item{--} if $\nu=1$ then $\tilde f$ is tangent to the identity in all characteristic directions of~$f$;
in other points of the exceptional divisor the eigenvalues of the differential of~$\tilde f$ are all equal to~$1$ but the differential is not diagonalizable.
\smallskip
\noindent This means that we can always repeat the previous construction blowing-up $\tilde f$ at a characteristic direction of~$f$; this will be important in the sequel.

As a first application of the blow-up construction, let us use it for describing the dynamics of dicritical maps.
If $f\in\End_1(\C^n,O)$ is dicritical, Theorem~\EcalleHakim\ yields a parabolic curve tangent to all directions outside a hypersurface of~$\P^{n-1}(\C)$ (notice that all directors are zero),
and the same holds for $f^{-1}$. One can then summarize the situation as follows:

\newthm Proposition \BMuno: ([Br1, 2]) Let $f\in\End_1(\C^n,O)$ be a dicritical germ tangent to the identity  of order~$\nu\ge 1$. 
Write $P_{\nu+1}(z)=p(z)z$, and let $D=\{[v]\in\P^{n-1}(\C)\mid p(v)=0\}$. Then there are two open sets $U^+$,~$U^-\subset\C^n\setminus\{O\}$ such that:
\smallskip
\itm{(i)} $\bar{U^+\cup U^-}$ is a neighborhood of $\P^{n-1}(\C)\setminus D$ in the blow-up $\widetilde{\C^n}$ of~$O$; 
\itm{(ii)} the orbit of any $z\in U^+$ converges to the origin tangent to a direction $[v]\in \P^{n-1}(\C)\setminus D$;
\itm{(iii)} the inverse orbit (that is, the orbit under~$f^{-1}$) of any $z\in U^-$ converges to the origin tangent to a\break\indent direction $[v]\in \P^{n-1}(\C)\setminus D$.

Coming back to the general situation, when $f\in\End_1(\C^n,O)$ is tangent to the identity its blow-up $\tilde f$ fixes pointwise the exceptional divisor; more precisely, the fixed point set of~$\tilde f$ is the full transform of the fixed point set of~$f$, and in particular $\tilde f$ has a at least a hypersurface of fixed points. 
This is a situation important enough to deserve a special notation.

\newdef Let $E$ be a connected (possibly singular) hypersurface in a complex manifold~$M$. We shall denote by $\End(M,E)$ the set of germs about~$E$ of holomorphic self-maps of~$M$ fixing pointwise~$E$.

If $E$ is a hypersurface in a complex manifold~$M$, we shall denote by $\ca O_M$ the sheaf of holomorphic functions on~$M$, and by $\ca I_E$ the subsheaf of functions vanishing on~$E$. Given $f\in\End(M,E)$, $f\not\equiv \id_M$, take $p\in E$. For every $h\in\ca O_{M,p}$, the germ $h\circ f$ is well-defined,
and $h\circ f-h\in\ca I_{E,p}$. Following [ABT1] (see also [ABT2, 3]), we can then introduce a couple of important notions. 

\newdef Let $E$ be a connected hypersurface in a complex manifold~$M$. Given $f\in\End(M,E)$, $p\in E$ and $h\in\ca O_{M,p}$, let $\nu_f(h)=\max\{\mu\in\N\mid
h\circ f-h\in\ca I_{E,p}^\mu\}$. Then the {\sl order of contact} $\nu_f$ of~$f$ with the identity along~$E$ is
$$
\nu_f= \min\{\nu_f(h)\mid h\in\ca O_{M,p}\}\;;
$$  
it can be shown ([ABT1]) that $\nu_f$ does not depend on~$p\in E$. Furthermore, we say that $f$ is {\sl tangential} if $\min\{\nu_f(h)\mid h\in\ca I_{E,p}\}>\nu_f$ for some (and hence any; see again [ABT1]) $p\in E$.

Let $(z_1,\ldots,z_n)$ be local coordinates in~$M$ centered at~$p\in E$, and $\ell\in\ca I_{E,p}$ a reduced equation of~$E$ at~$p$ (that is, a generator of~$\ca I_{E,p}$). If $(f_1,
\ldots,f_n)$ is the expression of $f$ in local coordinates, it turns out [ABT1] that we can write
$$
f_j(z)=z_j+\ell(z)^{\nu_f}g_j^o(z)
\neweq\eqtdue
$$
for $j=1,\ldots, n$, where there is a $j_0$ such that $\ell$ does not divide $g_{j_0}^o$; furthermore, $f$ is tangential if and only if $\nu_f(\ell)>\nu_f$. 

\newrem If $E$ is smooth at~$p$, we can choose local coordinates so that locally $E$ is given by $\{z_1=0\}$, that is $\ell=z_1$. Then we can write
$$
f_j(z)=z_j+z_1^{\nu_f}g_j^o(z)
$$
with $z_1$ not dividing some $g_j^o$; and $f$ is tangential if $z_1$ divides $g_1^o$, that is if
$f_1(z)=z_1+z_1^{\nu_f+1}h_1^o(z)$. More generally, if $E$ has a normal crossing at~$p$ with $1\le r\le n$ smooth branches, then we can choose local coordinates so that $\ell=z_1\cdots z_r$, so that $f_j(z)=z_j+(z_1\cdots z_r)^{\nu_f}g_j^o(z)$ with some $g_{j_0}^o$ not divisible by $z_1\cdots z_r$; in this case $f$ is tangential if and only if $z_j$ divides $g_j^o$ for $j=1,\ldots, r$. In particular, in the terminology of [A2] $f$ is tangential if and only if it is nondegenerate and $b_f=1$.

\newdef We say that $p\in E$ is a {\sl singular point} for~$f\in\End(M,E)$ (with respect to~$E$) if $g_1^o(p)=\cdots=g_n^o(p)$ in~\eqtdue; it turns out [ABT1] that this definition is independent of the local coordinates. Furthermore, the {\sl pure order} (or {\sl pure multiplicity})~$\nu_o(f,E)$ of~$f$ along~$E$ at~$p$ is
$$
\nu_o(f,E)=\min\{\ord_O(g_1^o),\ldots,\ord_O(g_n^o)\}\;.
$$
It is easy to see that the pure order does not depend on the local coordinates;
in particular, $p$ is singular for $f$ with respect to~$E$ if and only if $\nu_o(f,E)\ge 1$. If $E$ is the fixed point set of~$f$ at~$p$ then we shall talk of the {\sl pure order}~$\nu_o(f)$ of~$f$ at~$p$.

\newrem When $f$ is the blow-up of a germ $f_o\in\End_1(\C^n,O)$ tangent to the identity of order~$\nu\ge 1$, then \eqtuno\ implies that:
\smallskip
\item{--} $f$ is tangential if and only if $f_o$ is not dicritical; in particular, in this case being tangential is a generic condition;
\item{--} $\nu_f=\nu$ if $f_o$ is not dicritical, and $\nu_f=\nu+1$ if $f_o$ is dicritical;
\item{--} if $f_o$ is non dicritical, then $[v]\in\P^{n-1}(\C)$ is singular for~$f$ if and only if it is a characteristic direction of~$f_o$. 

Using the notion of singular points we can generalize Proposition~\Hakimz\ as follows:

\newthm Proposition \HakimA: ([ABT1]) Let $E\subset M$ be a hypersurface in a complex manifold~$M$, and $f\in\End(M,E)$, $f\not\equiv\id_M$, tangential to~$E$. Let $p\in E$ be a smooth point of~$E$ which is non-singular for~$f$. Then no infinite orbit of $f$ can stay arbitrarily close to $p$, that is, there exists a neighborhood~$U$ of~$p$ such that for all $q\in U$ either the orbit of~$q$ lands on~$E$ or $f^{n_0}(q)\notin U$ for some $n_0\in\N$. In particular, no infinite orbit is converging to~$p$.

More generally, we have:

\newthm Proposition \singint: ([AT1]) Let $f\in\End(\C^n,O)$ be of the form
$$
f_j(z)=\cases{z_j+z_j\left(\prod_{h=1}^r z_h^{\nu_h}\right)g_j(z)& for $1\le j\le r$,\cr
\noalign{\smallskip}
z_j+\left(\prod_{h=1}^r z_h^{\nu_h}\right)g_j(z)& for $r+1\le j\le n$,\cr}
\neweq\eqsinguno
$$
for suitable $1\le r<n$, with $\nu_1,\ldots,\nu_r\ge 1$ and $g_1,\ldots,g_n\in\ca O_{\C^n,O}$. Assume that
$g_{j_0}(O)\ne 0$ for some $r+1\le j_0\le n$. Then no infinite orbit can stay arbitrarily
close to $O$.


A very easy example of this phenomenon, promised at the end of the previous section, is the following:

\newex Let $f(z,w)=(z,w+z^2)$. Then $f$ is tangent to the identity at the origin; the fixed point set is $\{z=0\}$, and thus $O$ is not an isolated fixed point. We have $f^k(z,w)=(z,w+kz^2)$; therefore
all orbits outside the fixed point set escape to infinity, and in particular no orbit converges to the origin. Notice that this germ has only one characteristic direction, which is degenerate (and tangent to the fixed point set). Moreover, $f$ is tangential with order of contact~2 to its fixed point set, but the origin is not singular.

After these generalities, in the rest of this section we specialize to the case $n=2$ and to tangential maps (because of Remark~3.2; see anyway [ABT1] for information on the dynamics of non-tangential maps). Take $f\in\End(M,E)$, where $M$ is a complex surface and $E\subset M$ is a 1-dimensional curve smooth at~$p\in E$, and assume that $f$ is tangential to~$E$ with order of contact~$\nu_f\ge 1$. Then we can choose local coordinates centered at~$p$ so that we can write
$$
\cases{
f_1(z)=z_1+z_1^{\nu_f+1}h_1^o(z)\;,\cr
f_2(z)=z_2+z_1^{\nu_f}g_2^o(z)\;,\cr}
\neweq\eqtqua
$$
where $z_1$ does not divide $g_2^o$; notice that $h^o_1(0,\cdot)={\de g^o_1\over\de z_1}(0,\cdot)$, where $g^o_1=z_1h^o_1$. In particular, $O$ is singular if and only if $g_2^o(O)=0$. We then introduce the following definitions:

\newdef Let $f\in\End(M,E)$ be written in the form \eqtqua. Then:
\smallskip
\item{--} the {\sl multiplicity}~$\mu_p$ of~$f$ along~$E$ at~$p$  is $\mu_p=\ord_0\bigl(g_2^o(0,\cdot)\bigr)$, so that $p$ is a singular point if and only if $\mu_p\ge 1$;
\item{--} the {\sl transversal multiplicity}~$\tau_p$ of~$f$ along~$E$ at~$p$ is $\tau_p=\ord_0\bigl(h_1^o(0,\cdot)\bigr)$; 
\item{--} $p$ is an {\sl apparent singularity} if $1\le\mu_p\le\tau_p$;
\item{--} $p$ is a {\sl Fuchsian singularity} if $\mu_p=\tau_p+1$;
\item{--} $p$ is an {\sl irregular singularity} if $\mu_p>\tau_p+1$;
\item{--} $p$ is a {\sl non-degenerate singularity} if $\mu_p\ge 1$ but $\tau_p=0$;
\item{--} $p$ is a {\sl degenerate singularity} if $\mu_p$, $\tau_p\ge 1$;
\item{--} the {\sl index}~$\iota_p(f,E)$ of $f$ at~$p$ along~$E$ is
$$
\iota_p(f,E)=\nu_f\hbox{Res}_0 {h^o_1(0,\cdot)\over g^o_2(0,\cdot)}\;;
$$
\item{--} the {\sl induced residue} $\hbox{Res}^0_p(f)$ of $f$ along~$E$ at~$p$ is
$$
\hbox{Res}^0_p(f)=-\iota_p(f,E)-\mu_p\;.
$$
\smallskip
\noindent It is possible to prove (see [A2, ABT1, AT3]; notice that our index is $\nu_f$ times the residual index introduced in~[A2]) that these definitions are independent of the local coordinates; see also Remark~4.3.

\newrem Recalling \eqtuno, we see that if $f$ is obtained as the blow-up of a non-dicritical map~$f_o$, and $E$ is the exceptional divisor of the blow-up, then:
\smallskip
\item{--} the multiplicity of $[v]$ as characteristic direction of~$f_o$ is equal to the multiplicity of $f$ along~$E$ at $[v]$; 
\item{--} $[v]$ is a degenerate/non-degenerate characteristic direction of $f_o$ if and only if it is a degenerate/non-degenerate singularity of~$f$.
\smallskip
\noindent Furthermore, if we write $P_{\nu+1,1}(1,w)=\sum_{k=0}^{\nu+1} a_k w^k$ and $P_{\nu+1,2}(1,w)=\sum_{k=0}^{\nu+1} b_k w^k$ then $[1:0]$ is a characteristic direction if and only if
$b_0=0$, non-degenerate if and only if moreover $a_0\ne 0$, and (setting $b_{\nu+2}=0$)
$$
{h_1^o(0,\zeta)\over g_2^o(0,\zeta)}={1\over\zeta}{a_0+\sum_{k=1}^{\nu+1}a_k\zeta^k
\over (b_1-a_0)+\sum_{k=1}^{\nu+1}(b_{k+1}-a_k)\zeta^k}\;.
$$
So if $b_1\ne a_0$ we have 
$$
\mu_O=1,\quad 
\iota_{[1:0]}(f,E)={\nu a_0\over b_1-a_0}\;,\quad \Res^0_{[1:0]}(f)={(\nu-1)a_0+b_1\over a_0-b_1}\;;
$$
moreover, if $a_0\ne 0$ then $[1:0]$ is a non-degenerate characteristic direction with director $\alpha=(b_1-a_0)/\nu a_0$. More generally, we have $\tau_{[1:0]}=\ord_0\bigl(P_{\nu+1,1}(1,w)\bigr)$, $\mu_{[1:0]}=\ord_0\bigl(P_{\nu+1,2}(1,w)-wP_{\nu+1,1}(1,w)\bigr)$ and
$$
\iota_{[1:0]}={\nu a_{\mu-1}\over b_\mu-a_{\mu-1}}\;,\quad \Res^0_{[1:0]}={(\nu-\mu)a_{\mu-1}+\mu b_\mu\over a_{\mu-1}-b_\mu}\;,
$$
where $\mu=\mu_{[1:0]}$. In particular we obtain:
\smallskip
\item{--} if $[v]$ is non-degenerate characteristic direction of~$f_o$ with director~$\alpha\ne 0$ then 
$$
\iota_{[v]}(f,E)={1\over\alpha}\;; 
$$
\item{--} $[v]$ is a non-degenerate characteristic direction with non-zero director for~$f_o$ if and only if it is a Fuchsian singularity of multiplicity~1 for $f$.

Residues and indices are important for two reasons. First of all, we have the following {\sl index theorem:}

\newthm Theorem \index: (Abate, 2001 [A2]; Abate-Bracci-Tovena, 2004 [ABT1]) Let $E\subset M$ be a smooth compact Riemann surface in a complex surface~$M$. Let $f\in\End(M,E)$, $f\not\equiv \id_M$, be tangential with order of contact~$\nu$; denote by $\Sing(f)\subset E$ the finite set of singular points of~$f$ in~$E$. Then
$$
\sum_{p\in\Sing(f)}\iota_p(f)=\nu c_1(N_E)\;,\quad
\sum_{p\in\Sing(f)}\Res^0_p(f)=-\chi(E)\;,
$$
where $c_1(N_E)$ is the first Chern class of the normal bundle $N_E$ of $E$ in~$M$, and $\chi(E)$ is the Euler characteristic of~$E$. In particular, when $f$ is the blow-up of a nondicritical germ tangent to the identity and $E=\P^1(\C)$ is the exceptional divisor we have
$$
\sum_{p\in\Sing(f)}\iota_p(f)=-\nu \;,\quad
\sum_{p\in\Sing(f)}\Res^0_p(f)=-2\;.
$$

\newrem Bracci and Tovena [BT] have defined a notion of index at non-smooth points of~$E$ 
allowing the generalization of Theorem~\index\ to non necessarily smooth compact Riemann surfaces, where in the statement $c_1(N_E)$ is replaced by the self-intersection~$E\cdot E$.

The second reason is that the index can be used to detect the presence of parabolic curves. To state this precisely, we need a definition.

\newdef Let $f\in\End_1(\C^2,O)$ be tangent to the identity. We say that $O$ is a {\sl corner} if the germ of the fixed point set at the origin is locally reducible, that is has more than one irreducible component.

Then

\newthm Theorem \FAz: (Abate, 2001 [A2]) Let $E\subset M$ be a smooth Riemann surface in a complex surface~$M$, and take $f\in\End(M,E)$ tangential. Let $p\in E$ be a singular point, not a corner, such that $\iota_p(f)\notin \Q^+\cup\{0\}$. Then there exists a Fatou flower with~$\nu_f$ petals for $f$ at~$p$.

\newthm Corollary \FA: ([A2]) Let $f\in\End_1(\C^2,O)$ be tangent to the identity, and assume that
$O$ is a nondicritical singular point. Let $[v]\in\P^1(\C)$ be a characteristic direction, and $\tilde f$ the blow-up of~$f$. If $[v]$ is not a corner for~$\tilde f$ and $\iota_{[v]}(\tilde f)\notin \Q^+\cup\{0\}$ then there exists a Fatou flower for $f$ tangent to~$[v]$.

Theorem~\Abate\ is then a consequence of Theorem~\index\ and Corollary~\FA. Indeed, take $f\in\End_1(\C^2,O)$ tangent to the identity with an isolated fixed point at the origin. If $O$ is dicritical, we can directly apply Theorem~\EcalleHakim. Assume then $O$ non-dicritical, and let $\tilde f\in\End(\widetilde{\C^2},E)$ be the blow-up of~$f$. Since $O$ is non-dicritical, $\tilde f$ is tangential; Theorem~\index\ then implies that at least one characteristic direction~$[v]$ has negative index. Since $O$ is an isolated fixed point, the fixed point set of $\tilde f$ coincides with the exceptional divisor; therefore $[v]$ is not a corner, and Corollary~\FA\ yields the Fatou flower we were looking for. 

\newrem Bracci and Degli Innocenti (see [B, D]), using the definition of index introduced in [BT], have shown that Theorem~\FAz\ still holds when $E$ is not smooth at~$p$. Bracci and Suwa [BS] have also obtained a version of Theorem~\FAz\ when $M$ has a (sufficiently tame) singularity at~$p$. 

\newex Let $f(z,w)=(z+z^2w+zw^2+w^4,w+zw^2+z^4)$. Then $f$ is tangent to the identity at the origin of order~2, and the origin is an isolated fixed point. Furthermore, $f$ is non-dicritical and it has (see Example~2.2) two characteristic directions, $[v_1]=[1:0]$ and $[v_2]=[0;1]$, both degenerate. Working as in Remark~3.4 it is easy to see that $[v_1]$ is an irregular singularity of multiplicity~3 with index~$-2$ and induced residue~$-1$, and that $[v_2]$ is an apparent singularity of multiplicity~1, vanishing index, and induced residue~$-1$. In particular, $f$ admits a Fatou flower with 2 petals tangent to~$[v_1]$. 

\newex Let $f(z,w)=(z+w^2,w+z^3)$. Then $f$ is tangent to the identity at the origin of order~1, and the origin is an isolated fixed point. Furthermore, $f$ is non-dicritical with only one characteristic direction $[v]=[1:0]$, which is degenerate of multiplicity~3, Fuchsian, with index~$-1$ and induced residue~$-2$. Therefore $f$ admits a Fatou flower with one petal tangent to~$[v]$; compare with Example~2.3.

There are still instances where Theorem~\FAz\ cannot be applied:

\newex Let $f(z,w)=(z+zw+w^3,w+2w^2+bz^3)$ with $b\ne 0$. This map is tangent to the identity, with order 1, and the origin is an isolated fixed point. Moreover, it has two characteristic directions: $[1:0]$, degenerate Fuchsian with multiplicity~2, index~$1$ and induced residue~$-3$; and $[0:1]$, non-degenerate Fuchsian with multiplicity~1, index~$-2$ and induced residue~1. Theorem~\EcalleHakim\ (as well as Corollary~\FA) yields a Fatou flower tangent to~$[0:1]$; on the other hand, none of the results proven up to now say anything about direction $[1:0]$. 

However, a deep result by Molino gives the existence of a Fatou flower in the latter case too:

\newthm Theorem \Molinoz: (Molino, 2009 [Mo]) Let $E\subset M$ be a smooth Riemann surface in a complex surface~$M$, and take $f\in\End(M,E)$ tangential with order of contact~$\nu$. Let $p\in E$ be a singular point, not a corner, such that $\nu_o(f)=1$ and $\iota_p(f)\ne 0$. Then there exists a Fatou flower for $f$ at~$p$. More precisely:
\smallskip
\itm{(i)} if $p$ is an irregular singularity, or a Fuchsian singularity with $\iota_p(f)\ne \nu\mu_p$, then there exists a Fatou flower for $f$ with $\nu+\tau_p(\nu+1)$ petals;
\itm{(ii)} if $p$ is a Fuchsian singularity with $\iota_p(f)=\nu\mu_p$ then there exists a Fatou flower for $f$ with $\nu$ petals. 

\newrem Even more precisely, when $p$ is Fuchsian with $\mu_p\ge 2$ and $\iota_p(f)=\nu\mu_p$ then Molino constructs parabolic curves defined on the connected components of a set of the form
$$
D_{\nu+1-{1\over\mu_p},\delta}=\bigl\{\zeta\in\C\bigm| |\zeta^{r+1-1/\mu_p}(\log\zeta)^{1-1/\mu_p}
-\delta|<\delta\bigr\}\;,
$$
which has at least $\nu$ connected components with the origin in the boundary. 

\newthm Corollary \Molino: ([Mo]) Let $f\in\End(\C^2,O)$ be tangent to the identity of order~$\nu\ge 1$, and assume that $O$ is a nondicritical singular point. Let $[v]\in\P^1(\C)$ be a characteristic direction, and $\tilde f$ the blow-up of~$f$. If $[v]$ is not a corner for~$\tilde f$, $\nu_o(\tilde f)=1$ and $\iota_{[v]}(\tilde f)\ne 0$ then there exists a Fatou flower for $f$ with at least $\nu$ petals tangent to~$[v]$.

The assumption on the pure order in these statements seems to be purely technical; so it is natural to advance the following

\newthm Conjecture \conj: Let $E\subset M$ be a smooth Riemann surface in a complex surface~$M$, and take $f\in\End(M,E)$ tangential of order of contact~$\nu$. Let $p\in E$ be a singular point, not a corner, such that $\iota_p(f)\ne 0$. Then there exists a Fatou flower for $f$ at~$p$. 

See Section~5, and in particular (5.3), for examples of systems having Fatou flowers at singular points with vanishing index.
 
Instrumental in the proofs of Theorems~\FAz\ and~\Molinoz\ is a reduction of singularities statement. We shall need a few definitions:

\newdef Let $f\in\End_1(\C^n,O)$ be tangent to the identity. A {\sl modification} of
 $f$ is a $\tilde f\in\End(M,E)$ obtained as the lifting of~$f$ to a finite sequence of blow-ups, where the first one is centered in~$O$ and the remaining ones are centered in singular points of the intermediate lifted maps contained in the exceptional
 divisor. A modification is {\sl non-dicritical} if none of the centers of the blow-ups is dicritical. Associated to a modification $\tilde f\in\End(M,E)$ of~$f$ we have a holomorphic map
 $\pi\colon M\to\C^n$ such that $\pi^{-1}(O)=E$, $\pi|_{M\setminus E}$ is a biholomorphism
 between $M\setminus E$ and $\C^n\setminus\{O\}$, and $f\circ\pi=\pi\circ\tilde f$. The
 exceptional divisor $E$ is the union of a finite number of copies of $\P^{n-1}(\C)$, crossing
 transversally. 

\newdef Let $f\in\End_1(M,p)$ be tangent to the identity, where $M$ is a complex surface. In local coordinates centered at~$p$, we can write $f(z)=z+\ell(z)g^o(z)$, where $\ell=\gcd(f_1-z_1,f_2-z_2)$ is defined up to units. We shall say that $p$ is an {\sl irreducible singularity} if:
\smallskip
\item{(a)} $\ord_p(\ell)\ge 1$ and $\nu_o(f)=1$; and
\item{(b)} if $\lambda_1$, $\lambda_2$ are the eigenvalues of the linear part of~$g^o$ then either
\itemitem{$(\star_1)$} $\lambda_1$, $\lambda_2\ne 0$ and $\lambda_1/\lambda_2$, $\lambda_2/\lambda_1\notin\N$, or
\itemitem{$(\star_2)$} $\lambda_1\ne 0$, $\lambda_2=0$.

It turns out that there always exists a modification with only dicritical or irreducible singularities:

\newthm Theorem \redsing: (Abate, 2001 [A2]) Let $f\in\End_1(\C^2,O)$ be tangent to the identity, and assume that $O$ is a singular point. Then there exists a non-dicritical modification $\tilde f\in\End(M,E)$ of~$f$ such that the singular points of~$\tilde f$ on~$E$ are either irreducible or dicritical.

\newdef Let $f\in\End_1(\C^2,O)$ be tangent to the identity. The modification of $f$ satisfying the conclusion of Theorem~\redsing\ obtained with the minimum number of blow-ups is the {\sl minimal resolution} of~$f$.

It is easy to see that the techniques of the proof of Theorem~\EcalleHakim\ yield the existence of a Fatou flower at dicritical singularities, and at irreducible singularities of type $(\star_1)$ which are not a corner; then the proof of Theorem~\FAz\ amounts to showing that if the index is not a non-negative rational number then the minimal resolution contains at least a singularity which is either dicritical or of type $(\star_1)$ and not a corner. The proof of Theorem~\Molinoz.(i) consists in showing that, under those hypotheses, the minimal resolution must contain a non-degenerate singularity, which is not a corner and where one can apply Theorem~\EcalleHakim; the proof of Theorem~\Molinoz.(ii) requires instead a technically hard extension of Theorem~\EcalleHakim.

See also [Ro4, 7] and [LS] for other approaches to resolution of singularities for germs tangent to the identity in arbitrary dimension, and [AT2], [AR] for the somewhat related problem of the identification of formal normal forms for germs tangent to the identity. 

\smallsect 4. Parabolic domains

Theorem~\Hakimb\ yields conditions ensuring the existence of parabolic domains attached to a non-degenerate characteristic direction. In dimension~2, Vivas has found conditions ensuring the existence of a parabolic domain attached to Fuchsian and irregular degenerate characteristic directions, and Rong has found conditions ensuring the existence of a parabolic domain attached to apparent degenerate characteristic directions. Very recently, Lapan [L2] has extended Rong's approach to cover more types of degenerate characteristic directions. 

More precisely, Vivas has proved the following result: 

\newthm Theorem \Vivasdue: (Vivas, 2012 [V1]) Let $f\in\End_1(\C^2,O)$ be tangent to the identity of order~$\nu\ge 1$, with $O$ nondicritical. 
Let $[v]\in\P^1(\C)$ be a degenerate characteristic direction, and $\tilde f$ the blow-up of~$f$. Denote by $\mu\ge 1$ the multiplicity , by $\tau\ge 0$ the transversal multiplicity, by $\iota\in\C$ the index, and by $\nu_o\ge 1$ the pure order of~$\tilde f$ at~$[v]$. 
Assume that either
\smallskip
\itm{(a)} $[v]$ is Fuchsian (thus necessarily $\tau\ge 1$ because $[v]$ is degenerate) and 
$$
\Re\iota+\tau>0\;,\qquad \left|\iota +{\tau\over2}-{\nu\mu\over2}\right|>{\tau\over 2}+{\nu\mu\over 2}\;;
$$
or
\itm{(b)} $[v]$ is Fuchsian, $\nu_o=1$ and
$$
\left|\iota-{\mu\nu\over2}\right|<{\mu\nu\over 2}\;;
$$ 
or
\itm{(c)} $[v]$ is Fuchsian, $\nu_o>1$ and
$$
\Re\iota+\tau>0\;, \qquad \left|\iota-{(\nu+1)\tau\over 2}\right|>{(\nu+1)\tau\over 2}\;;
$$
or
\itm{(d)} $[v]$ is irregular. 
\smallskip\noindent Then there is a parabolic domain attached to~$[v]$. 

See also Remark~6.4 for a comment about the conditions on~$\iota$ and~$\tau$.

To state Rong's theorem, consider a germ  $f\in\End_1(\C^2,O)$ tangent to the identity of order~$\nu\ge 1$,
and assume that $[1:0]$ is a characteristic direction of~$f$. Then we can write 
$$
\cases{ f_1(z,w)=z+az^{r+1}+O(z^{r+2})+w\alpha(z,w),\cr 
f_2(z,w)=w+bz^\nu w+dz^{s+1}+O(z^{s+2})+O(wz^{\nu+1})+w^2\beta(z,w),\cr}
\neweq\eqquno
$$
with $\nu\le r\le+\infty$, $\nu+1\le s\le+\infty$, $\ord_O(\alpha)\ge\nu$, $\ord_O(\beta)\ge\nu-1$, and $a\ne 0$ if $r<+\infty$ (respectively, $d\ne 0$ if $s<+\infty$). The characteristic direction $[1:0]$ is non-degenerate if and only if $r=\nu$; in this case the director is given by ${1\over\nu}\left({b\over a}-1\right)$. On the other hand, saying that $[1:0]$ is degenerate with $b\ne 0$ is equivalent to saying that $r>\nu$ and that $[1:0]$ has multiplicity~1 and transversal multiplicity at least~1; in particular, in this case it is an apparent singularity. Then Rong's theorem can be stated as follows:

\newthm Theorem \Rongsette: (Rong, 2014 [Ro9]) Let $f\in\End_1(\C^2,O)$ be tangent to the identity of order~$\nu\ge 1$ and written in the form~$\eqquno$, so that $[1:0]$ is a characteristic direction. Assume that $r>\nu$ and $b\ne 0$, so that $[1:0]$ is an apparent degenerate characteristic direction. Suppose furthermore that $s>r$, and that $b^2/a\notin\R^+$ if $r=2\nu$. Then there is a parabolic domain attached to~$[1:0]$. 

To state Lapan's result we need to introduce a few definitions. 

\newdef Let $f\in\End_1(\C^2,O)$ be tangent to the identity of order~$\nu\ge 1$ with homogeneous expansion \equzero. We say that $[v]\in\P^1(\C)$ is a {\sl characteristic direction of degree}~$s\ge\nu+1$ if it is a characteristic direction of $P_{\nu+1},\ldots,P_s$. We shall say that it is {\sl non-degenerate in degree}~$r+1$, with $\nu<r<s$, if it is degenerate for $P_{\nu+1},\ldots, P_r$ and non-degenerate for~$P_{r+1}$.

For instance, if $f$ is in the form \eqquno\ with $s<+\infty$, then $[1:0]$ is a characteristic direction of degree~$s$. If furthermore $r+1\le s$ then it is non-degenerate in degree~$r+1$.

\newdef Let $f\in\End_1(\C^2,O)$ be tangent to the identity of order~$\nu\ge 1$ with homogeneous expansion \equzero. Assume that $[1:0]\in\P^1(\C)$ is a characteristic direction of degree $s\ge\nu+1$. Given $\nu+1\le j\le s$, the {\sl $j$-order} of~$[1:0]$ is the order of vanishing at~$0$ of $P_{j,2}(1,\cdot)$, where $P_{j}=(P_{j,1},P_{j,2})$. We say that $[1;0]$ is {\sl of order one in degree~$t+1$}, with $\nu\le t<s$, if the $j$-order of~$[1:0]$ is larger than one for~$\nu+1\le j\le t$ and of $(t+1)$-order exactly equal to~$1$.

For instance, if $b\ne 0$ in \eqquno\ then $[1:0]$ is of order one in degree~$\nu+1$. More generally, if $b=0$ and $[1:0]$ is of order one in degree~$t+1$ then we can replace the term
$O(wz^{\nu+1})$ by $O(wz^{t+1})$.

Assume that $[1:0]$ is a characteristic direction of degree~$s<+\infty$, non-degenerate in degree~$r+1\le s$ and of order one in degree~$t+1\le s$. Then we can write
$$
\cases{f_1(z,w)=z+az^{r+1}+O(z^{r+2})+w\alpha(z,w),\cr
f_2(z,w)=w+bz^{t+1} w+dz^{s+1}+O(z^{s+2})+O(wz^{t+2})+w^2\beta(z,w),\cr}
\neweq\eqqdue
$$
with $abd\ne 0$. Then Lapan's theorem can be stated as follows:

\newthm Theorem \Lapan: (Lapan, 2015 [L2]) Let $f\in\End_1(\C^2,O)$ be tangent to the identity of order~$\nu\ge 1$ and written in the form~$\eqqdue$, so that $[1:0]$ is a characteristic direction of degree $s<+\infty$, non-degenerate in degree $\nu+1\le r+1\le s$, and of order one in degree~$t+1\le s$. Assume that $t\le r$ and $s>r+t-\nu$. Suppose furthermore that either
\smallskip
\itm{(i)} $r\ne t$, $2t$, or
\itm{(ii)} $r=t$ and $\Re(b/a)>0$, or
\itm{(iii)} $r=2t$ and $b^2/a\notin\R^+$.
\smallskip
\noindent Then there is a parabolic domain attached to~$[1:0]$. 

The assumptions of Theorem~\Rongsette\ imply that $[1:0]$ is a characteristic direction of degree $s$, non-degenerate in degree $\nu+1<r+1\le s$, and of order one in degree~$t+1=\nu+1$; therefore Theorem~\Rongsette\ is a particular case of Theorem~\Lapan.

Parabolic domains are often used to build Fatou-Bieberbach domains, that is proper subsets of~$\C^n$ biholomorphic to~$\C^n$; see, e.g., [V2], [SV] and references therein.

\smallsect 5. The formal infinitesimal generator 

A different approach to the study of parabolic curves in~$\C^2$ has been suggested by Brochero-Mart\'\i nez, Cano and L\'opez-Hernanz [BCL],
and further developed by C\^amara and Sc\'ardua [CaS] and by L\'opez-Hernanz and S\'anchez [LS]. It consists in using the formal infinitesimal generator of a germ tangent to the identity. To describe their approach, we need to introduce several definitions.

\newdef We shall denote by $\widehat{\ca O}_n=\C[\![ z_1,\ldots,z_n]\!]$ the space of formal power series in~$n$ variables. The {\sl order}~$\ord(\hat\ell)$ of $\hat\ell\in\widehat{\ca O}_n$ is the 
lowest degree of a non-vanishing term in the Taylor expansion of~$\hat\ell$. A {\sl formal map} is 
a $n$-tuple of formal power series in $n$ variables; the space of formal maps will be denoted by $\widehat{\ca O_n}^n$. We shall 
denote by $\widehat{\End}(\C^n,O)$ the set of formal maps with vanishing constant term; by $\widehat{\End}_1(\C^n,O)$ the subset of formal maps tangent to the identity,
and by $\widehat{\End}_\nu(\C^n,O)$ the subset of formal maps tangent to the identity of order at least~$\nu\ge 1$. 

\newdef We shall denote by $\ca X_n$ the space of germs at the origin of holomorphic vector fields in~$\C^n$. A {\sl formal vector field} is an expression of the form  
$\hat X=\hat X_1{\de\over\de z_1}+\cdots+\hat X_n{\de\over\de z_n}$ where $\hat X_1,\ldots,\hat X_n\in\widehat{\ca O}_n$ are the {\sl components} of~$\hat X$. The space of {\sl formal vector fields} will be denoted by $\widehat{\ca X}_n$.  The {\sl order} $\ord(\hat X)$ of $\hat X\in\widehat{\ca X}_n$ is the minimum among the orders of its components. 
We put $\widehat{\ca X}^k_n=\{\hat X\in\widehat{\ca X}_n\mid \ord(\hat X)\ge k\}$. If $\hat X\in\widehat{\ca X}^k_n$, the {\sl principal part} of~$\hat X$ will be the unique polynomial
homogeneous vector field $H_k$ of degree exactly~$k$ such that $\hat X-H_k\in\widehat{\ca X}_n^{k+1}$. A {\sl characteristic direction} for~$\hat X$ is an invariant line for~$H_k$.

\newrem There is a clear bijection between $\widehat{\ca X}_n$ and $\widehat{\ca O_n}^n$ obtained by associating to a formal vector field the $n$-tuple of its components; so we shall sometimes identify formal vector fields and formal maps without comments.
In particular, this bijection preserves characteristic directions.

If $X\in\ca X_n$ is a germ of holomorphic vector field vanishing at the origin (that is, of order at least~1), the associated time-1 map~$f_X$ will be a well defined germ in~$\End(\C^n,O)$, that can be recovered as follows (see, e.g., [BCL]):
$$
f_X=\sum_{k\ge 0}{1\over k!}X^{(k)}(\id)\;,
\neweq\eqqdue
$$
where $X^{(k)}$ is the $k$-th iteration of $X$ thought of as derivation of $\End(\C^n,O)$. Now, not every germ in $\End(\C^n,O)$ can be obtained as a time-1 map of a convergent vector field (see, e.g., [IY, Theorem~21.31]). However, it turns out that the right-hand side of \eqqdue\ is well-defined as a formal map for all $X\in\widehat{\ca X}_n^1$.

\newdef The {\sl exponential map} $\exp\colon\widehat{\ca X}_n^1\to\widehat{\End}(\C^n,O)$ is defined by the right-hand side of \eqqdue.

When $k\ge 2$, if $\hat X\in\widehat{\ca X}_n^k$ has principal part~$H_k$ then it is easy to check that 
$$
\exp(\hat X)=\id+H_k+\hbox{h.o.t.}\;;
\neweq\eqqumez
$$
in particular, the exponential of a formal vector field of order~$k$ is a formal map tangent to the identity of order~$k-1$. Takens (see, e.g., [IY, Theorem~3.17]) has shown that on the formal
level the exponential map is bijective:

\newthm Proposition \Takens: The exponential map $\exp\colon\widehat{\ca X}_n^{\nu+1}\to\widehat{\End}_\nu(\C^n,O)$ is bijective for all $\nu\ge 1$.

\newdef If $\hat f\in\widehat{\End}_\nu(\C^n,O)$, the unique formal vector field $\hat X\in\widehat{\ca X}_n^{\nu+1}$ such that $\exp(\hat X)=\hat f$ is the {\sl formal infinitesimal generator} of~$\hat f$.

The idea now is to read properties of a holomorphic germ tangent to the identity from properties of its formal infinitesimal generator, using Theorem~\EcalleHakim\ as bridge for going back from the formal side to the holomorphic side.

The first observation is that if $\pi\colon(\widetilde{\C^2},E)\to(\C^2,O)$ is the blow-up of the origin, $\hat X\in\widehat{\ca X}_2^2$ is a formal vector field and $[v]\in E$ is
a characteristic direction of (the principal part of)~$\hat X$ then we can find a formal vector field $\hat X_{[v]}\in\widehat{\ca X}_2^2$ such that $d\pi(\hat X_{[v]})=\hat X\circ\pi$.
This lifting is compatible with the exponential in the following sense:

\newthm Proposition \BCLu: ([BCL]) Let $f\in\End_1(\C^2,O)$ be tangent to the identity with formal infinitesimal generator $\hat X\in\widehat{\ca X}_2^2$, and let $\tilde f\in\End(\widetilde{\C^2},E)$ be the blow-up of~$f$. Let $[v]\in E$ be a characteristic direction of~$f$, and denote by $\tilde f_{[v]}$ the germ of~$\tilde f$ at~$[v]$. Then
$\tilde f_{[v]}=\exp(\hat X_{[v]})$.

In particular, Brochero-Mart\'\i nez, Cano and L\'opez-Hernanz's proofs of Theorems~\Abate\ and~\FAz\ go as follows: let $\hat X\in\widehat{\ca X}_2^2$ be the formal infinitesimal generator of $f\in\End_1(\C^2,O)$ with an isolated fixed point (so that $\hat X$ has an isolated singular point at the origin). Then the formal version of Camacho-Sad's theorem [CS] (see also [Ca]) shows that we can find a finite composition $\pi\colon(M,E)\to(\C^2,O)$ of blow-ups at singular points and a smooth point $p\in E$ such that the lifting $\hat X_p$ of~$\hat X$, in suitable coordinates centered at~$p$ adapted to~$E$ (in the sense that $E$ is given by the equation $\{z=0\}$ near~$p$), has the expression
$$
\hat X_p(z,w)=z^m\left(\bigl(\lambda_1 z+O(z^2)\bigr){\de\over\de z}+\bigl(\lambda_2 w+O(z)\bigr){\de\over\de w}\right) 
$$
with $\lambda_1\ne 0$, $\lambda_2/\lambda_1\notin\Q^+$ and $m\ge\ord(\hat X)-1$. Then $\exp(\hat X_p)$ has the form
$$
\exp(\hat X_p)(z,w)=\bigl(z+\lambda_1 z^{m+1}+O(z^{m+2}),w+\lambda_2 z^mw+O(z^{m+1})\bigr)\;,
$$
which has a non-degenerate characteristic direction transversal to~$E$ --- and hence a Fatou flower outside the exceptional divisor. By Proposition~\BCLu, $\exp(\hat X_p)$
is the blow-up of $\exp(\hat X)=f$; therefore projecting this Fatou flower down by~$\pi$ we get a Fatou flower for~$f$.

In [CaS] and [LS] this approach has been pushed further showing how to relate formal separatrices and parabolic curves. 

\newdef A {\sl formal curve}~$\hat C$ in~$(\C^2,0)$ is a reduced principal ideal of~$\widehat{\ca O}_2$. Any generator of the ideal is an {\sl equation} of the curve; the equation is defined up to an unit in~$\widehat{\ca O}_2$. The {\sl tangent cone} of a formal curve~$\hat C$ is the set of zeros of the homogeneous part of least degree of any equation of~$\hat C$; the {\sl tangent directions} to~$\hat C$ are the points in~$\P^1(\C)$ determined by the tangent cone.

\newrem It is known that a formal curve $\hat C$ is irreducible if and only if it has a unique tangent direction. 

\newdef Let $\hat X\in\widehat{\ca X}_2^2$. A {\sl singular formal curve} for~$\hat X$ is
a formal curve $\hat C=(\hat\ell)$ such that $\hat X=\hat\ell \hat X_1$ for some $\hat X\in\widehat{\ca X}_2^1$. 
A {\sl formal separatrix} of~$\hat X$ is a formal curve $\hat C=(\hat\ell)$ such that $\hat X(\hat\ell)\in(\hat\ell)$. Clearly singular formal curves are formal separatrices.  

The corresponding notions for germs tangent to the identity are:

\newdef Let $f\in\End_1(\C^2,O)$. A formal curve $\hat C=(\hat\ell)$ is a {\sl formal separatrix} for~$f$ if $\hat\ell\circ f\in(\hat\ell)$. In particular, this means that $f$ acts by composition on~$\widehat{\ca O}_2/(\hat\ell)$; if the action is the identity, we say that $\hat C$ is {\sl completely fixed} by~$f$. Notice that $\hat C$ is completely fixed by~$f$  if and only if we can write $f=\id+\hat\ell \hat g$ for some $\hat g\in\widehat{\ca O}_2^2$.  

\newthm Proposition \CaSuno: ([CaS]) Let $\hat X\in\widehat{\ca X}_2^2$ be the formal infinitesimal generator of
$f\in\End_1(\C^2,O)$. Then:
\smallskip
\itm{(i)} a formal curve is a formal separatrix for $f$ if and only if it is a formal separatrix for~$\hat X$;
\itm{(ii)} a formal curve is completely fixed for $f$ if and only if it is a singular formal curve for $\hat X$;
\itm{(iii)} a completely fixed curve for $f$ always has a convergent equation;
\itm{(iv)} the tangent directions to a formal separatrix are characteristic directions for~$f$, and the tangent directions\break\indent to a completely fixed curve are degenerate characteristic directions for~$f$.

Let $\hat C=(\hat\ell)$ be a formal curve, and $[v]\in\P^1(\C)$ a tangent direction to~$\hat C$. If $\pi\colon(\widetilde{\C^2},E)\to(\C^2,O)$ is the blow-up of the origin, 
we can find a formal curve $\pi^*\hat C_{[v]}=(\hat\ell_{[v]})$ at~$[v]$ such that $\hat\ell_{[v]}=\hat\ell\circ\pi$; the tangent directions to~$\pi^*\hat C_{[v]}$ are higher order tangent directions of~$\hat C$. This construction can be iterated, and it gives a way of lifting formal curves along a finite sequence of blow-ups. Using this idea, and a generalization of Hakim's technique, L\'opez-Hernanz and S\'anchez have been able to prove the following 

\newthm Theorem \LSuno: (L\'opez-Hernanz and S\'anchez, [LS]) Let $f\in\End_1(\C^2,O)$ be a germ tangent to the identity admitting a formal separatrix $\hat C$ not completely fixed. Then $f$ or $f^{-1}$ (or both) admit a parabolic curve tangent to (a tangent direction of) $\hat C$. 

\newrem Actually, [LS, Theorem~1] gives the more precise statement that the parabolic curve $\phe\colon D\to\C^2$ is {\sl asymptotic} to the formal separatrix $\hat C$. This means that there exists a formal parametrization $\hat\gamma\in\widehat{\ca O_1}^2$ of $\hat C$ such that for every $N\in\N$
there exists $c_N>0$ such that 
$$
|\phe(\zeta)-(J_N\hat\gamma)(\zeta)|\le c_N|\zeta|^{N+1}
$$
for all $\zeta\in D$, where $J_N\hat\gamma$ is the $N$-th jet of $\hat\gamma$. A formal parametrization of $\hat C$ is a formal map $\hat\gamma\in\widehat{\ca O_1}^2$ such that $g\in\hat C$ if and only if $g\circ\hat\gamma\equiv 0$. 

\newrem In [CaS] C\^amara and Sc\'ardua claimed that under the hypotheses of Theorem~\LSuno\ 
$f$ must admit a parabolic curve tangent to $\hat C$. Unfortunately, the core of their argument was [CaS, Proposition~2.12], and in its proof they forgot to consider vector fields of the form $\hat X^o(z,w)=z(1+\lambda w^p){\de\over\de z}+w^{p+1}{\de\over\de w}$, where their approach does not work.

The proof of Theorem~\LSuno\ has three steps. First of all, the authors show that, assuming the existence of a formal separatrix not completely fixed, after a finite number of blow-ups the germ $f$ can be brought in the following normal form:
$$
\cases{
f_1(z,w)=z+z^{\nu+p+1}\bigl(\lambda+\psi(z,w)\bigr)\;,\cr
f_2(z,w)=w+z^\nu\bigl(b(z)+a(z)w+O(w^2)\bigr)\;,\cr}
\neweq\eqcitre
$$
with $p\ge 0$, $\lambda\ne 0$, $\ord_O\psi\ge 1$, $a(0)\ne 0$ and $b(0)=0$. For a germ in this form, a {\sl real attracting direction} is a $\tau\in S^1$ such that $\tau^{\nu+p}\lambda=-1$. Then, generalizing Hakim's proof of Theorem~2.2, the authors show that if $f$ has a real attracting direction 
$\tau$ such that $\Re\left({a(0)\over\lambda\tau^p}\right)<0$ then $f$ admits a parabolic curve, that turns out to be asymptotic to $\hat C$. Finally, they show that at least one between $f$ and $f^{-1}$ have a real attracting direction satisfying the given condition.

Notice that a germ in the form \eqcitre\ has pure order~1, but vanishing index if $p\ge 1$; so we cannot apply Theorem~3.7. On the other hand, there are germs tangent to the identity admitting parabolic curves thanks to Theorem~3.7 but without formal separatrices not completely fixed:

\newex Let $f=\exp(z^\nu X^o)$, with $\nu\ge 2$ and 
$$
\hat X^o(z,w)=z\bigl(\lambda+ A(z,w)\bigr){\de\over\de z}+\bigl(z+\lambda w+ B(z,w)\bigr){\de\over\de w}\;,
$$
with $\lambda\ne 0$, $\ord_O(A)\ge 1$ and $\ord_O(B)\ge 2$. Then
$$
f(z,w)=\Bigl( z+z^{\nu+1}\bigl(\lambda+A(z,w)\bigr)+O(z^{2\nu+1}),w+z^\nu\bigr(z+\lambda w+B(z,w)\bigr)+O(z^{2\nu})\Bigr)\;.
$$
The germ $f$ has pure order~1, and the linear part of~$g^o$ is not diagonalizable; since the index of~$f$ at~$O$ along the fixed point set is~1, Theorem~\Molinoz\ yields a Fatou flower. Furthermore, $f$ has a unique (degenerate) characteristic direction, $[0:1]$. Blowing up and looking in the coordinates centered at~$[0:1]$ we get
$$
\tilde f(u,w)=\Bigl(u+u^{\nu+1} w^\nu\bigl(-u+w\hat A(u,w)\bigr), w+u^\nu w^{\nu+1}\bigl(\lambda+u+w\hat B(u,w)\bigr)\Bigr)=\exp(u^\nu w^\nu \tilde X^o)\;,
$$
where $\tilde X^o= u\hat C(u,w){\de\over\de u}+w\bigl(\lambda+\hat D(u,w)\bigr){\de\over\de w}$ with $\ord_O(\hat C)$,~$\ord_O(\hat D)\ge 1$. Since the linear part of~$\tilde X^o$ is diagonalizable, $X^o$ has exactly two formal separatrices, necessarily given by the axes. It follows that all formal separatrices of $\tilde f$ are completely fixed;
so to $\tilde f$ we cannot apply Theorem~\LSuno, but $\tilde f$ still has a Fatou flower because $f$ does.

The paper [LS] also indicates a way to adapt these techniques to more than two variables. However, it should be kept in mind that [AT1] contains examples in $\C^3$ of germs tangent to the identity without parabolic curves asymptotic to formal separatrices. 

\smallsect 6. Homogeneous vector fields and geodesics.

None of the results presented up to now (with the partial exception of Proposition~\BMuno) describe the dynamics in a full neighborhood of the fixed point, and so in this sense they are not a complete generalization of the Leau-Fatou flower theorem. As far as we know, up to now the only techniques able to give results in a full neighborhood are the ones introduced in [AT3], that we shall briefly describe now.

We have seen that every germ tangent to the identity can be realized as the time-1 map of a formal vector field of order at least~2; and that a lot of information can be deduced from the principal part of this vector field, principal part which is a {\sl homogeneous} vector field. Furthermore, in dimension 1 Camacho-Shcherbakov theorem (see [C, Sh]) says that every germ tangent to the identity is locally topologically conjugated to the time-1 map of a homogeneous vector field. So time-1 maps of homogeneous vector fields clearly are an important class of examples; and the insights we obtain from their study (and, more generally, from the study of the real dynamics of homogeneous vector fields) can shed light on the dynamics of more general germs tangent to the identity.

The work described in [AT3] had exactly the aim of studying the real dynamics of homogenous vector fields in~$\C^n$; for the sake of clarity, here we shall summarize only the results in~$\C^2$ only, referring to [AT3] for more general statements. 

Let $H\in{\ca X}_2^{\nu+1}$ be a homogeneous vector field in~$\C^2$ of degree~$\nu+1\ge 2$. It clearly determines a homogeneous self-map of~$\C^2$ of the same degree;
in particular, we can adapt to~$H$ all the definitions we introduced for homogeneous self-maps (degenerate/non-degenerate characteristic directions, multiplicities, index, induced residue, being dicritical).

\newdef Let $H\in{\ca X}_2$ be a homogeneous vector field in~$\C^2$. A {\sl characteristic line} for~$H$ is a line $L_v=\C v$ which is $H$-invariant, that is such that $[v]\in\P^1(\C)$ is a characteristic direction.

If $L_v=\C v$ is a characteristic line then integral lines of~$H$ issuing from points in~$L_v$ stay inside~$L_v$. If $[v]$ is degenerate, $H$ vanishes identically along~$L_v$, and so the dynamics there is trivial. If $[v]$ is non-degenerate, then the dynamics inside~$L_v$ is one-dimensional, and can be summarized as follows:

\newthm Lemma \ATdcl: Let $[v]\in\P^1(\C)$ be a non-degenerate characteristic direction of a
homogeneous vector field $H=H_1{\de\over\de z_1}+H_2{\de\over\de z_2}\in{\ca X}_2^{\nu+1}$ of degree $\nu+1\ge 2$. Choose a representative $v\in\C^2$ so that $H(v)=v$. Then the real integral curve of $H$ issuing from $\zeta_0 v\in L_v$
is given by
$$
\gamma_{\zeta_0 v}(t)={\zeta_0 v\over(1-\zeta_0^\nu \nu t)^{1/\nu}}\;.
$$
In particular no (non-constant) integral curve is recurrent, and we have:
\smallskip
\item{\rm (a)} if $\zeta_0^\nu\notin\R^+$ then $\lim\limits_{t\to+\infty}
\gamma_{\zeta_0 v}(t)=O$;
\item{\rm (b)} if $\zeta_0^\nu\in\R^+$ then $\lim\limits_{t\to(\zeta_0^\nu\nu)^{-1}}\|\gamma_{\zeta_0 v}(t)\|=+\infty$.

Lemma \ATdcl\ completely describes the dynamics of dicritical homogeneous vector fields. More precisely, we see that every non-degenerate characteristic line contains a Fatou flower; thus in this case Theorem~\EcalleHakim\ becomes trivial, and we can shift our interest to the understanding of the dynamics outside the characteristic lines. To do so we need to introduce a new ingredient:

\newdef Let $\nabla^o$ be a meromorphic connection on $\P^1(\C)$ (see [IY] for an introduction to meromorphic connections), and denote by $\Sing(\nabla^o)$ the set of its poles. A {\sl geodesic} for~$\nabla^o$ is a smooth real curve $\sigma\colon I\to\P^1(\C)\setminus\Sing(\nabla^o)$ such that 
$$
\nabla^o_{\dot\sigma}\dot\sigma\equiv O\;.
$$ 

The main result allowing the understanding of the real dynamics of homogeneous vector fields is the following:

\newthm Theorem \ATctreb: (Abate-Tovena, 2011 [AT3]) Let $H\in{\ca X}_2^{\nu+1}$ be a non-dicritical homogeneous vector field of degree~$\nu+1\ge 2$ in~$\C^2$, and denote by $V_H$ the complement in~$\C^2$ of the characteristic lines of~$H$. Then there exists a meromorphic connection~$\nabla^o$ on~$\P^1(\C)$, whose poles are a (possibly proper) subset of the characteristic directions of~$H$, such that:
\smallskip
\item{\rm(i)} if
$\gamma\colon I\to V_H$ is a real integral curve of~$H$ then its
direction $[\gamma]\colon I\to\P^1(\C)$ is a geodesic for $\nabla^o$; conversely, 
\item{\rm(ii)} if $\sigma\colon I\to\P^1(\C)$ is a geodesic for  $\nabla^o$ then
there exists exactly $\nu$ real integral curves $\gamma_1,\ldots,\gamma_\nu\colon I\to V_H$\break\indent
of~$H$, differing only by the multiplication by a $\nu$-th root of unity, whose direction
is given by~$\sigma$, that is\break\indent such that $\sigma=[\gamma_j]$.

\newrem If $H=H_1{\de\over\de z_1}+H_2{\de\over\de z_2}\in{\ca X}_2^{\nu+1}$ is a homogeneous vector field of degree~$\nu+1\ge 2$, the meromorphic 1-form representing~$\nabla^o$ in the standard chart centered at~$0\in\P^1(\C)$ is given by ([AT3]) 
$$
\eta^o=-\left[\nu{H_1(1,\zeta)\over R(\zeta)}+{R'(\zeta)\over R(\zeta)}\right]\,d\zeta\;,
$$
where $R(\zeta)=H_2(1,\zeta)-\zeta H_1(1,\zeta)$; a similar formula, exchanging the r\^ole of $z_1$ and $z_2$, holds in the standard chart centered at~$\infty\in\P^1(\C)$.
In particular recalling \eqqumez, \eqtuno\ and Definition~3.7 we see that the poles of $\nabla^o$ are singular points for the blow-up $\tilde f$ of the time-1 map of~$H$,
and that the residue of $\nabla^o$ at a pole $p\in\P^1(\C)$ coincides with the induced residue of~$\tilde f$ at~$p$.\hfill\break\indent
Furthermore, in [AT3] we introduced another meromorphic connection~$\nabla$ defined on the $\nu$-th tensor power~$N_E^{\otimes\nu}$ of the normal bundle~$N_E$ of the exceptional divisor $E=\P^1(\C)$ in the blow-up of the origin in~$\C^2$. The meromorphic 1-form representing $\nabla$ in the standard chart centered at~$0\in\P^1(\C)$ is given by
$$
\eta=-\nu{H_1(1,\zeta)\over R(\zeta)}\,d\zeta\;.
$$
Therefore the poles of $\nabla$ are exactly the Fuchsian and irregular characteristic directions of~$\tilde f$, and the residue of~$\nabla$ at a pole~$p\in\P^1(\C)$ coincides with the opposite of the index of~$\tilde f$ at~$p$.

So the study of the real integral curves of~$H$ is reduced to the study of the geodesics of a meromorphic connection on~$\P^1(\C)$. This study is subdivided in two parts: the study of the global behavior of geodesics, and the study of the local behavior nearby the poles. It turns out that the global behavior is related to the induced residues, while the 
local behavior is mainly related to the index. To state our results we need a couple of definitions.

\newdef A geodesic $\sigma\colon[0,\ell]\to\P^1(\C)$ for a meromorphic connection $\nabla^o$ is {\sl closed}
if $\sigma(\ell)=\sigma(0)$ and $\sigma'(\ell)$ is a positive 
multiple of $\sigma'(0)$; it is {\sl periodic} if $\sigma(\ell)=\sigma(0)$
and $\sigma'(\ell)=\sigma'(0)$.

Contrarily to the case of Riemannian geodesics, closed geodesics are not necessarily periodic; see [AT3]. The (induced) residue allows to recognize closed and periodic geodesics:

\newthm Proposition \ATqunomez: ([AT3]) Let $\nabla^o$ be a meromorphic connection on 
$\P^1(\C)$, with poles $\{p_0,p_1,\ldots,p_r\}$, and set
$S=\P^1(\C)\setminus\{p_0,\ldots,p_r\}\subseteq\C$. Let $\sigma\colon
[0,\ell]\to S$ be a geodesic with $\sigma(0)=\sigma(\ell)$ and no
other self-intersections; in particular, $\sigma$ is an oriented
Jordan curve. Let $\{p_1,\ldots,p_g\}$ be the poles of $\nabla^o$
contained in the interior of~$\sigma$. Then $\sigma$ is a closed geodesic if and only if
$$
\sum_{j=1}^g\Re\Res_{p_j}(\nabla^o)=-1\;,
$$
and it is a periodic geodesic if and only if
$$
\sum_{j=1}^g\Res_{p_j}(\nabla^o)=-1\;.
$$
If $\sigma$ is closed, it can be extended to an
infinite length geodesic $\sigma\colon J\to S$, where $J$ is a half-line
(possibly $J=\R$). Moreover,
\smallskip
\itm{(i)}if $\sum\limits_{j=1}^g \Im\Res_{p_j}(\nabla)<0$ then 
$\sigma'(t)\to O$ as $t\to+\infty$ and $|\sigma'(t)|\to+\infty$ as 
$t$ tends to the other end of~$J$;
\itm{(ii)}if $\sum\limits_{j=1}^g \Im\Res_{p_j}(\nabla)>0$ then
$\sigma'(t)\to O$ as $t\to-\infty$ and 
$|\sigma'(t)|\to+\infty$ as $t$ tends to the other end of~$J$.

\newthm Corollary \ATccin: Let $\gamma\colon\R\to\C^2\setminus\{O\}$ be a
non-constant periodic integral curve of a
homogeneous vector field $H$ of degree~$\nu+1\ge 2$. Then the
characteristic directions $[v_1],\ldots,[v_g]\in\P^1(\C)$ surrounded
by~$[\gamma]$ satisfy
$$
\sum_{j=1}^g \Res^0_{[v_j]}(H)=-1\;,
$$
where $\Res^0_{[v_j]}(H)$ denotes the induced residue at~$[v_j]$ of the blow-up of the time-1 map of~$H$.

Closed but not periodic geodesics correspond to integral curves converging to the origin on one side and escaping to infinity on the other side; the convergence to the origin occurs along a spiral, and thus the time-1 map has orbits converging to the origin without being tangent to any direction. This can actually happen; see [AT3] for an example. 

\newdef Let $\sigma\colon I\to S$ be a curve in $S=\P^1(\C)\setminus\{p_0,
\ldots,
p_r\}$. A {\sl simple loop in~$\sigma$} is the restriction 
of $\sigma$ to a closed interval $[t_0,t_1]\subseteq I$ such that 
$\sigma|_{[t_0,t_1]}$ is a simple loop~$\tau$. If $p_1,\ldots, p_g$
are the poles of~$\nabla$ contained in the interior of~$\tau$,
we shall say that $\tau$ {\sl surrounds}~$p_1,\ldots,p_g$.

\newdef A {\sl saddle connection} for a meromorphic connection $\nabla^o$ 
on $\P^1(\C)$ is a maximal geodesic
$\sigma\colon(\eps_-,\eps_+)\to \P^1(\C)$
(with $\eps_-\in[-\infty,0)$ and $\eps_+\in(0,+\infty]$) such that
$\sigma(t)$ tends to a pole of~$\nabla^o$ both when $t\uparrow\eps_+$ and
when $t\downarrow\eps_-$. A {\sl graph of saddle connections} is a connected graph
in $\P^1(\C)$ made up of saddle connections.

Then we have a Poincar\'e-Bendixson type theorem, describing the asymptotic behavior of geodesics:

\newthm Theorem \ATqqua: (Abate-Tovena, 2011 [AT3]; Abate-Bianchi, 2014 [AB]) Let $\sigma\colon[0,\eps_0)\to S$ be a maximal
geodesic for a meromorphic connection~$\nabla^o$ on $\P^1(\C)$, where
$S=\P^1(\C)\setminus\{p_0,\ldots,p_r\}$ and $p_0,\ldots,p_r$ are the
poles of~$\nabla^o$. Then either
\smallskip
\itm{(i)} $\sigma(t)$ tends to a pole of~$\nabla^o$ as $t\to\eps_0$; or
\itm{(ii)} $\sigma$ is closed, and then surrounds poles $p_1,\ldots,p_g$
 with $\sum\limits_{j=1}^g\Re\Res_{p_j}(\nabla^o)=-1$; or
\item{\rm (iii)} the $\omega$-limit set of $\sigma$ in $\P^1(\C)$ is given by the support of a closed geodesic surrounding poles
$p_1,\ldots,p_g$ with $\sum\limits_{j=1}^g\Re\Res_{p_j}(\nabla^o)=-1$; or
\item{\rm (iv)} the $\omega$-limit set of $\sigma$ in $\P^1(\C)$ is a graph of saddle connections whose complement in~$\P^1(\C)$ has a connected component containing $p_1,\ldots,p_g$ with $\sum\limits_{j=1}^g\Re\Res_{p_j}(\nabla^o)=-1$; or
\itm{(v)} $\sigma$ intersects itself infinitely many times, and in this
case every simple loop of $\sigma$ surrounds a set of poles whose
sum of residues has real part belonging to $(-3/2,-1)\cup(-1,-1/2)$.
\smallskip
\noindent In particular, a recurrent geodesic either intersects itself
infinitely many times or is closed.

\newthm Corollary \cqua: Let $H$ be a homogeneous holomorphic vector
field on $\C^2$ of degree~$\nu+1\ge 2$, and let $\gamma\colon[0,\eps_0)\to\C^2$
be a recurrent maximal integral curve of $Q$. Then $\gamma$ is periodic or 
$[\gamma]\colon[0,\eps_0)\to\P^1(\C)$ intersects itself infinitely many times.

\newrem We have examples (see [AT3]) of cases (i), (ii), (iii) and (v), but not yet of case (iv). 

It is worthwhile to notice that the maximal geodesics of generic meromorphic connections will
behave as in case (i), because the other cases require that a particular relationships between the (induced) residues should hold. In particular, this means that the direction of a maximal real integral curve of a generic homogeneous vector field will go from a characteristic line to a characteristic line (possibly the same); the next step then consists in understanding what happens nearby characteristic lines. It turns out that we can find holomorphic normal forms for apparent and Fuchsian singularities, and formal normal forms for irregular singularities; and we shall see that the local behavior is mostly related to the index.

To key behind this local study is the following

\newthm Theorem \ATtre: (Abate-Tovena, 2011 [AT3]) Let $N_E$ be the normal bundle of the exceptional divisor of the blow-up~$(M,E)$ of the origin in~$\C^2$. Then for every $\nu\ge 1$ there exists a holomorphic $\nu$-to-$1$ covering $\chi_\nu\colon\C^2\setminus\{O\}\to N_E^{\otimes\nu}\setminus E$ satisfying $\pi\circ\chi_\nu(v)=[v]$, where $\pi\colon N_E^{\otimes\nu}\to E=\P^1(\C)$ is the canonical projection, such that for every homogeneous vector field
$H\in{\ca X}_2^{\nu+1}$ of degree~$\nu+1$ the push-forward $d\chi_\nu(H)$ defines a global holomorphic vector field $G$ on the total space of~$N_E^{\otimes\nu}$. In particular, a real curve $\gamma\colon I\to\C^2\setminus\{O\}$ is an integral curve for~$H$ if and only if $\chi_\nu\circ\gamma$ is an integral curve for~$G$.

\newdef The field $G$ is the {\sl geodesic field} associated to the homogeneous vector field~$H$. The reason of the name is that the projections on~$\P^1(\C)$ of the integral curves of~$G$ are geodesics for the connection~$\nabla^o$ associated to~$H$. 

The point is that the field $G$ has a form well suited to reduction to normal form. Indeed, if we denote by $\zeta$ the usual coordinate on~$\C\subset\P^1(\C)$ centered at the origin, and by $v$ the corresponding coordinate on the fibers of~$N_E^{\otimes\nu}$, which over~$\C$ is canonically trivialized, we have 
$$
G(\zeta,v)=R(\zeta)v{\de\over\de\zeta}+\nu H_1(1,\zeta) v^2{\de\over\de v}\;,
$$
where $R(\zeta)=H_2(1,\zeta)-\zeta H_1(1,\zeta)$ as before; and a similar formula holds in the usual coordinates centered at~$\infty\in\P^1(\C)$. In particular, we can read the multiplicity and the transversal multiplicity (and hence the type of singularity) of a characteristic direction of~$H$ in the order of vanishing of the components of the geodesic
field. 

Since $G$ is a vector field on the total space of a line bundle, it is natural to consider only changes of coordinates preserving the bundle structure, that is changes of coordinates of the form
$$
(\zeta,v)\mapsto \bigl(\psi(\zeta),\xi(\zeta)v\bigr)\;,
$$
where $\psi$ a germ of biholomorphism and $\xi$ is a non-vanishing holomorphic function. It turns out that these changes of coordinates are enough to obtain normal forms around apparent and Fuchsian singularities. 

For apparent singularities we have the following theorem:

\newthm Theorem \ATsuno: (Abate-Tovena, 2011 [AT3]) Let $[v]\in\P^1(\C)$ be an apparent characteristic direction of multiplicity~$\mu\ge 1$ of a homogeneous vector field $H\in{\ca X}_2^{\nu+1}$. Then we can find local coordinates centered at~$[v]$ such that in these coordinates the geodesic field~$G$ associated to~$H$ is given by
$$
G=\cases{\displaystyle
\zeta v {\de\over\de\zeta}&if $\mu=1$,\cr
\displaystyle \zeta^\mu(1+a \zeta^{\mu-1})v{\de\over\de\zeta}&for some $a\in\C$ if $\mu>1$.\cr}
$$
Furthermore, if $\mu>1$ then $a\in\C$ is a holomorphic invariant, the {\rm apparent index.}

In particular, around an apparent singularity the geodesic field $G$ is explicitly integrable. Studying the integral lines of~$G$ and rephrasing the results in terms of the integral curves of~$H$ we obtain the following corollary:

\newthm Corollary \ATstre: ([AT3]) Let $H\in{\ca X}_2^{\nu+1}$ be a homogeneous vector
field on $\C^2$ of degree~$\nu+1\ge 2$. 
Let $[v]\in\P^1(\C)$ be an apparent singularity 
of~$H$ of multiplicity~$\mu\ge 1$ (and apparent index $a\in\C$ if $\mu>1$). Then:
\smallskip
\itm{(i)} if the direction $[\gamma(t)]\in\P^1(\C)$ of an integral curve
$\gamma\colon[0,\eps)\to\C^2\setminus\{O\}$ of $H$ tends to $[v]$
as $t\to\eps$ then $\gamma(t)$ tends to a non-zero point
of the characteristic leaf $L_{v}\subset\C^2$;
\itm{(ii)} no integral curve of $H$ tends to the origin tangent to~$[v]$;
\itm{(iii)} there is an
open set of initial conditions whose integral curves tend to a non-zero 
point of~$L_{v}$;
\itm{(iv)} if $\mu=1$ or $\mu>1$ and $a\ne 0$ then $H$ admits periodic orbits of
arbitrarily long periods accumulating at\break\indent the origin.

In particular, in case (iv) the time-1 map of~$H$ has both periodic orbits accumulating at the origin (small cycles), when the period of the integral curve is rational, and orbits whose closure a is a closed Jordan curve, when the period of the integral curve is irrational; both phenomena cannot happen in one variable.

The holomorphic classification of Fuchsian singularities is the following: 

\newthm Theorem \ATsttre: (Abate-Tovena, 2011 [AT3]) Let $[v]\in\P^1(\C)$ be a Fuchsian characteristic direction of multiplicity~$\mu\ge 1$, transversal multiplicity~$\tau=\mu-1\ge 0$ and index~$\iota\in\C^*$ of a homogeneous vector field $H\in{\ca X}_2^{\nu+1}$. Then we can find local coordinates centered at~$[v]$ such that in these coordinates the geodesic field~$G$ associated to~$H$ is given:
\smallskip
\itm{(i)} if $\tau+\iota\notin\N^*$ by
$$
\zeta^{\mu-1}\left(\zeta v{\de\over\de\zeta}+\iota v^2{\de\over\de v}\right)\;;
$$
\itm{(ii)} if $n=\tau+\iota\in\N^*$ by
$$
\zeta^{\mu-1}\left(\zeta v{\de\over\de\zeta}+\iota v^2(1+a \zeta^n){\de\over\de v}\right)
$$
for a suitable $a\in\C$ which is a holomorphic invariant, the {\rm resonant index.}

When the resonant index is zero the integral curves of the geodesic field can be expressed in terms of elementary functions and easily studied. This is not the case when the resonant index is different from zero; however we are able to obtain the following description of the integral curves of~$H$ nearby Fuchsian characteristic directions:

\newthm Corollary \ATstcin: ([AT3])  Let $H\in{\ca X}_2^{\nu+1}$ be a homogeneous vector
field on $\C^2$ of degree~$\nu+1\ge 2$. 
Let $[v]\in\P^1(\C)$ be a Fuchsian singularity  
of~$H$ of multiplicity~$\mu\ge 1$, transversal multiplicity $\tau=\mu-1\ge 0$ and index $\iota\in\C^*$ (and resonant index $a\in\C$ if 
$\tau+\iota\in\N^*$). Then:
\smallskip
\itm{(i)} if the direction $[\gamma(t)]\in\P^1(\C)$ of an integral curve
$\gamma\colon[0,\eps)\to\C^2\setminus\{O\}$ of $H$ tends to $[v]$
as $t\to\eps$ and $\gamma$ is not contained in the characteristic leaf $L_{v}$ then 
\itemitem{\rm(a)} if $\tau+\Re\iota>0$ and $\left|\iota+{\tau\over 2}\right|>{\tau\over 2}$ then $\gamma(t)$ tends to the origin;
\itemitem{\rm(b)} if $\tau+\iota=0$, or $\tau+\Re\iota<0$, or $\tau+\Re\iota>0$ and $\left|\iota+{\tau\over 2}\right|<{\tau\over 2}$, then $\|\gamma(t)\|$ tends to $+\infty$;
\itemitem{\rm(c)} if $\tau+\Re\iota>0$ and $\left|\iota+{\tau\over 2}\right|={\tau\over 2}$ then $\gamma(t)$ accumulates a
circumference in $L_{v}$.
\smallskip
\noindent Furthermore there is
a neighbourhood $U\subset\P^1(\C)$ of $[v]$ such that an integral curve
$\gamma$ issuing from a point $z_0\in\C^2\setminus L_{v}$ with 
$[z_0]\in U\setminus\{[v]\}$ 
must have one of the following behaviors, where $\hat U=\{z\in\C^2
\setminus\{O\}\mid[z]\in U\}$:
\smallskip
\itm{(ii)} if $\tau+\Re\iota<0$ then
\smallskip
\itemitem{\rm (a)}either $\gamma(t)$ escapes $\hat U$, and this happens
for a Zariski open dense set of initial conditions in~$\hat U$; or
\itemitem{\rm(b)} $[\gamma(t)]\to[v]$ but $\|\gamma(t)\|\to+\infty$;
\smallskip
in particular, no integral curve outside $L_{v}$ converge to the origin tangent to~$[v]$;
\itm{(iii)} if $\tau+\Re\iota=0$ but $\tau+\iota\ne0$ then 
\smallskip
\itemitem{\rm(a)} either $\gamma(t)$ escapes $\hat U$; or
\itemitem{\rm(b)} $\gamma(t)\to O$ without being tangent to any direction,
and $[\gamma(t)]$ is a closed curve or accumulates a closed
curve in~$\P^1(\C)$ surrounding $[v]$; or
\itemitem{\rm(c)} $\|\gamma(t)\|\to+\infty$ without being tangent to any direction,
and $[\gamma(t)]$ is a closed curve in~$\P^1(\C)$ surrounding $[v]$;
\smallskip
in particular, no integral curve outside $L_{v}$ converge to the origin tangent to~$[v]$;
\itm{(iv)} if $\tau+\iota=0$ then 
\smallskip
\itemitem{\rm(a)} either $\gamma(t)$ escapes $\hat U$, and this 
happens for an open set $\hat U_1\subset\hat U$ of initial conditions; or
\itemitem{\rm(b)} $[\gamma(t)]\to[v]$ with $\|\gamma(t)\|\to+\infty$,
and this happens for an open set $\hat U_2\subset\hat U$ of initial conditions
such that $\hat U_1\cup\hat U_2$ is dense in $\hat U$; or
\itemitem{\rm(c)} $\gamma$ is a periodic integral curve with $[\gamma]$
surrounding $[v]$;
\smallskip
in particular, no integral curve outside $L_{v}$ converge to the origin tangent to~$[v]$, but we have periodic\break\indent integral curves of arbitrarily long period accumulating the origin; 
\itm{(v)} if $\tau+\Re\iota>0$ and $a=0$ then $[\gamma(t)]\to[v]$
for an open 
dense set~$\hat U_0$ of initial conditions in~$\hat U$, and $\gamma$
escapes $\hat U$ for $z\in\hat U\setminus\hat U_0$; moreover,
\smallskip
\itemitem{\rm (a)} if $\left|\iota+{\tau\over 2}\right|>{\tau\over 2}$ then $\gamma(t)\to O$
tangent to~$[v]$ for all $z\in\hat U_0$;
\itemitem{\rm (b)} if $\left|\iota+{\tau\over 2}\right|<{\tau\over 2}$ then $\|\gamma(t)\|\to
+\infty$ tangent to $[v]$ for all $z\in\hat U_0$;
\itemitem{\rm(c)} if $\left|\iota+{\tau\over 2}\right|={\tau\over 2}$ then $\gamma(t)$
accumulates a circumference in $L_{v}$.

\newrem We conjecture that Corollary~\ATstcin.(v) should hold also when $a\ne 0$.

\newrem This result must be compared with Theorems~\Hakimb\ and \Vivasdue. We already noticed that a non-degenerate characteristic direction $[v]$ with non-zero
director~$\delta$ is a Fuchsian singularity of multiplicity~1, and hence transversal multiplicity~0. Then 
Corollary~\ATstcin\ says that if $\Re\iota<0$ (that is $\Re\delta<0$)
then no orbit of the time-1 map of $H$ outside of~$L_{v}$
converges to the origin tangent to~$[v]$,
whereas if $\Re\iota>0$ (that is $\Re\delta>0$) and $\iota$ is not a positive integer (or $a=0$ if $\iota\in\N^*$) 
then the orbits under the time 1-map of $H$ converge to the origin tangent
to~$[v]$ for an open (and dense in a conical neighbourhood of~$[v]$)
set of initial conditions, providing the existence of a parabolic domain in accord with Theorem~\Hakimb.\hfill\break\indent
If instead $\tau>0$ and $\tau+\iota\notin\N^*$ (or $a=0$ if $\tau+\iota\in\N^*$) then Corollary~\ATstcin\ yields a parabolic basin when $\tau+\Re\iota>0$ and $\left|\iota+{\tau\over 2}\right|>{\tau\over 2}$, which is a condition strictly weaker than the condition found in Theorem~\Vivasdue.(a); this suggests that there might be room for improvement in
the statement of the latter theorem.

Putting all of this together we can finally have a completely description of the dynamics for a substantial class of examples. For instance, we get the following:

\newthm Corollary \ATstsei: ([AT3]) Let $H\in{\ca X}_2^{\nu+1}$ be a 
homogeneous vector field on $\C^2$ of degree~$\nu+1\ge 2$. Assume
that $H$ is non-dicritical and all its characteristic directions are Fuchsian of multiplicity~$1$, 
Assume moreover that for no set of characteristic
directions the real part of the sum of the induced residues belongs to the interval $(-3/2,-1/2)$.
Let $\gamma\colon[0,\eps_0)\to\C^2$ be a maximal integral curve of $H$. Then:
\smallskip
\itm{(a)} If $\gamma(0)$ belongs to a characteristic leaf~$L_{v_0}$, then the image of $\gamma$ is contained 
in~$L_{v_0}$. Moreover, either $\gamma(t)\to O$ (and this happens for
a Zariski open dense set of initial conditions in~$L_{v_0}$), or $\|\gamma(t)\|\to+\infty$.
\itm{(b)} If $\gamma(0)$ does not belong to a characteristic leaf then either
\itemitem{\rm (i)} $\gamma$ converges to the origin tangentially to a characteristic 
direction $[v_0]$ whose index has positive real part; or
\itemitem{\rm (ii)} $\|\gamma(t)\|\to+\infty$ tangentially to a characteristic
direction $[v_0]$ whose index has negative real part. 
\smallskip
\noindent Furthermore, case {\rm (i)} happens for a Zariski open set of initial conditions.

\newrem The conditions in Corollary~\ATstsei\ imply that there must be at least one index with positive real part. Indeed, if the multiplicity is~1 then the induced residue is one less the opposite of the index. So assuming that no sum of $1\le g\le\nu+2$ induced residues has real part belonging to the interval $(-3/2,-1/2)$ is equivalent to saying that no sum of $1\le g\le\nu+2$ indices has real part belonging to the interval $({1\over 2}-g,{3\over 2}-g)$. Assume, by contradiction, that no index has positive real part; than the real part of all of them should be less than~$-1/2$. So the real part of the sum of two indices must be less than~$-1<-1/2$; so it should be less than $-3/2$. Arguing by induction on~$g$ one then shows that the sum of the real part of all indices should be less than ${1\over 2}-(\nu+2)=-\nu-3/2<-\nu$, against Theorem~\index, contradiction.

\newex Corollary~\ATstsei\ describes completely the dynamics of most vector fields of the form
$$
H(z,w)=\bigl(\rho z^2+(1+\tau)zw\bigr){\de\over\de z}+\bigl((1+\rho)zw+\tau w^2\bigr){\de\over\de w}\;.
$$
Indeed such a vector field has exactly three Fuchsian characteristic directions with multiplicity~1 and indices respectively $\rho$, $\tau$ and $-1-\rho-\tau$; so the conditions required by Corollary~\ATstsei\ are satisfied as soon as $\Re\rho$, $\Re\tau\notin(-1/2,1/2)$ and $\Re(\rho+\tau)\notin(-3/2,-1/2)$.

\smallsect 7. Other systems with parabolic behavior

Another situation where Fatou flowers can exist is when the eigenvalues of the differential are all equal to~1 but the differential is not necessarily diagonalizable. The reason is that we can reduce to the tangent to the identity case by using a suitable sequence of blow-ups:

\newthm Theorem \Azuno: (Abate, 2000 [A1])
Let $f\in\End(\C^n,O)$ be such that all eigenvalues of $df_O$ are equal to~$1$. Then there
exist a complex $n$-dimensional manifold~$M$, a holomorphic projection $\pi\colon
M\to\C^n$, a canonical point ${\bf e}\in M$ and a germ around $\pi^{-1}(O)$ of
holomorphic self-map $\tilde f\colon M\to M$ such that:
\smallskip
\item{\rm(i)}$\pi$ restricted to $M\setminus\pi^{-1}(O)$ is a biholomorphism between
$M\setminus\pi^{-1}(O)$ and $\C^n\setminus\{O\}$;
\item{\rm(ii)}$\pi\circ\tilde f=f\circ\pi$;
\item{\rm(iii)}$\bf e$ is a fixed point of~$\tilde f$ where $\tilde f$ is
tangent to the identity.

It should be remarked that the projection $\pi\colon M\to\C^n$ is obtained as a sequence of blow-ups whose centers are not necessarily reduced to points, and depend on the Jordan structure of~$df_O$. Furthermore $\pi$ is chosen in such a way that the interesting part of the dynamics of~$\tilde f$ is outside the exceptional divisor~$E$ (which is not in general pointwise fixed by~$\tilde f$), allowing the study of the dynamics of~$f$ by means of the dynamics of~$\tilde f$. For instance, we can get the following 

\newthm Proposition \Azdue: ([A1]) Let $f=(f_1,\ldots,f_n)\in\End(\C^n,O)$ be such that $df_O$ is not diagonalizable with all eigenvalues equal to~$1$. Without loss of generality we can suppose that $df_O$ is in Jordan form with $\rho$ blocks of order respectively $\mu_1\ge\mu_2\ge\cdots\ge\mu_\rho\ge 1$, where $\mu_1+\cdots+\mu_\rho=n$. Assume that $\mu_1>\mu_2$ and that the coefficient of $(z_1)^2$ in~$f_{\mu_1}$ is not zero. Then $f$ admits a parabolic curve tangent to~$[1:0:\cdots:0]$.

In dimension~2, using the tools introduced in [A2], one can get a cleaner result:

\newthm Corollary \Aztre: ([A2]) Let $f\in\End(\C^2,O)$ be such that $df_O$ is a Jordan block with eigenvalue~$1$, and assume that the origin is an isolated fixed point. Then $f$ admits a parabolic curve tangent to~$[1:0]$.

See [Ro6] (and [A3] for a particular example) for a detailed study of the existence of parabolic domains for germs in~$\End(\C^2,O)$ with non-diagonalizable differential. 

Finally, I would like to mention that parabolic curves, parabolic domains and Fatou flowers appear also in non-parabolic dynamical systems. This is not surprising 
in semi-parabolic systems, that is when the eigenvalues of the differential are either equal to~1 or in modulus strictly less than~1 (see, e.g., [N, H1, Ri2, U1, U2, Ro5]),
or in quasi-parabolic systems, where the eigenvalues of the differential are either equal to~1 or have modulus equal to~1 (see, e.g., [BM, Ro1, Ro2]). On the other hand,
a recent surprising discovery is that they also appear in {\sl multi-resonant systems,} whose differential is not parabolic at all but whose eigenvalues satisfies some resonance relation; see, e.g., [BZ, BRZ, RV, BR] for the main results of this very interesting theory. 

\bigskip\noindent

\setref{ABT3}
\beginsection References.

\art A1 M. Abate: Diagonalization of non-diagonalizable discrete holomorphic dynamical systems! Amer. J. Math.! 122 2000 757-781

\art A2 M. Abate: The residual index and the dynamics of holomorphic maps tangent to the identity! Duke Math. J.! 107 2001 173-206

\art A3 M. Abate: Basins of attraction in quadratic dynamical systems with a Jordan fixed point! Nonlinear Anal.! 51 2002 271-282

\coll A4 M. Abate: Discrete holomorphic local dynamical systems! Holomorphic dynamical systems! Eds. G. Gentili, J. Gu\'enot, G. Patrizio. Lect. Notes in Math. 1998, Springer, Berlin, 2010, pp. 1--55

\pre AB M. Abate, F. Bianchi: A Poincar\'e-Bendixson theorem for meromorphic connections on compact Riemann surfaces! Preprint, arXiv:1406.6944! 2014

\art ABT1 M. Abate, F. Bracci, F. Tovena: Index theorems for holomorphic
self-maps! Ann. of Math.! 159 2004 819-864

\art ABT2 M. Abate, F. Bracci, F. Tovena: Index theorems for holomorphic
maps and foliations! Indiana Univ. Math. J.! 57 2008 2999-3048

\art ABT3 M. Abate, F. Bracci, F. Tovena: Embeddings of submanifolds and
normal bundles! Adv. Math.! 220 2009 620-656

\art AR M. Abate, J. Raissy: Formal Poincar\'e-Dulac renormalization for holomorphic germs! Disc. Cont. Dyn. Syst.! 33 2013 1773-1807

\art AT1 M. Abate, F. Tovena: Parabolic curves in $\C^3$! Abstr. Appl. Anal.! 2003 2003 275-294

\art AT2 M. Abate, F. Tovena: Formal classification of holomorphic maps tangent to the identity! Disc. Cont. Dyn. Sys.! Suppl. 2005 1-10

\art AT3  M. Abate, F. Tovena: Poincar\'e-Bendixson theorems for meromorphic connections and holomorphic homogeneous vector fields! J. Diff. Eq.! 251 2011 2612-2684

\coll ArR M. Arizzi, J. Raissy: On \'Ecalle-Hakim's theorems in holomorphic dynamics! Progress in complex dynamics! A. Bonifant et al. eds., Princeton University Press, Princeton, 2014, pp. 387--449

\art B F. Bracci: The dynamics of holomorphic maps near curves of fixed points! Ann. Scuola Norm. Sup. Pisa! 2 2003 493-520

\art BM F. Bracci, L. Molino: The dynamics near quasi-parabolic fixed points of holomorphic diffeomorphisms in $\C^2$! Amer. J. Math.! 126 2004 671-686

\art BRZ F. Bracci, J. Raissy, D. Zaitsev: Dynamics of multi-resonant biholomorphisms! Int. Math. Res. Not. IMRN! {} 2013 4772-4797

\art BR F. Bracci, F. Rong: Dynamics of quasi-parabolic one-resonant biholomorphisms! J. Geom. Anal.! 24 2014 1497-1508

\art BS F. Bracci, T. Suwa: Residues for singular pairs and dynamics of biholomorphic maps of singular surfaces! Internat. J. Math.! 15 2004 443-466

\art BT F. Bracci, F. Tovena: Residual indices of holomorphic maps relative to singular curves of fixed points on surfaces! Math. Z.! 242 2002 481-490

\art BZ F. Bracci, D. Zaitsev: Dynamics of one-resonant biholomorphisms! J. Eur. Math. Soc.! 15 2013 179-200

\art Br1 F.E. Brochero-Mart\'\i nez: Groups of germs of analytic diffeomorphisms in $(\C^2,0)$!
J. Dyn. Control Sys.! 9 2003 1--32
 
\art Br2 F.E. Brochero-Mart\'\i nez: Din\'amica de difeomorfismos dicr\'\i ticos en $(\C^n,0)$! Rev. Semin. Iberoam. Mat.! 3 2008 33-40

\art BCL F.E. Brochero Mart\'\i nez, F. Cano, L. L\'opez-Hernanz: Parabolic curves for diffeomorphisms in $\C^2$!  Publ. Mat.! 52 2008 189-194

\art C C. Camacho: On the local structure of conformal mappings and holomorphic
vector fields! Ast\'e\-risque! 59--60 1978 83-94

\art CS C. Camacho, P.Sad: Invariant varieties through singularities of holomorphic vector fields! Ann. of Math! 115 1982 579-595

\art CaS L. C\^amara, B. Sc\'ardua: A Fatou type theorem for complex map germs! Conf. Geom. Dyn.! 16 2012 256-268

\art Ca J. Cano: Construction of invariant curves for singular holomorphic vector fields! Proc. Am. Math. Soc.! 125 1997 2649-2650

\art D F. Degli Innocenti: Holomorphic dynamics near germs of singular curves! Math. Z.! 251 2005 943-958

\book E  J. \'Ecalle: Les fonctions r\'esurgentes. Tome I\negthinspace
I\negthinspace I: L'\'equation du pont et la classification analytique des
objects locaux! Pr\'epublications Math. Orsay {\bf 85-05,} Universit\'e de Paris-Sud, Orsay,
1985 

\art F1 P. Fatou: Sur les \'equations fonctionnelles, I! Bull. Soc. Math. France!
47 1919 161-271

\art F2 P. Fatou: Sur les \'equations fonctionnelles, I\negthinspace I! Bull.
Soc. Math. France! 48 1920 33-94

\art F3 P. Fatou: Sur les \'equations fonctionnelles, I\negthinspace
I\negthinspace I! Bull. Soc. Math. France! 48 1920 208-314

\book GH P. Griffiths, J. Harris: Principles of algebraic geometry! Wiley, New York, 1978

\art H1 M. Hakim: Attracting domains for semi-attractive transformations of $\C^p$! Publ. Mat.! 38 1994 479-499

\art H2 M. Hakim: Analytic transformations of $(\C^p,0)$ tangent to the
identity! Duke Math. J.! 92 1998 403-428

\pre H3 M. Hakim: Transformations tangent to the identity. Stable pieces of manifolds! Preprint! 1997

\book IY Y. Ilyashenko, S. Yakovenko: Lectures on analytic differential equations! American Mathematical Society, Providence, 2008

\art L1 S. Lapan: Attracting domains of maps tangent to the identity whose only characteristic direction is non-degenerate! Int. J. Math.! {24 \rm 1350083} 2013 1-27

\pre L2 S. Lapan: Attracting domains of maps tangent to the identity in $\C^2$ with characteristic direction of multiple degrees! Preprint, arXiv:1501.00244! 2015

\art Le L. Leau: \'Etude sur les \'equations fonctionelles \`a une ou plusieurs variables!
Ann. Fac. Sci. Toulouse! 11 1897 E1-E110.

\pre LS L. L\'opez-Hernanz, F. Sanz S\'anchez: Parabolic curves of diffeomorphisms asymptotic to formal invariant curves! Preprint, arXiv:1411.2945! 2014

\book M  J. Milnor: Dynamics in one complex variable! Annals of Mathematics Studies, 160, Princeton University Press, Princeton, NJ, 2006

\art Mo L. Molino: The dynamics of maps tangent to the identity and with non vanishing index! Trans. Am. Math. Soc.! 361 2009 1597-1623

\art N Y. Nishimura: Automorphismes analytiques admettant des sous-vari\'et\'es de points fix\'es attractives dans la direction transversale! J. Math. Kyoto Univ.! 23 1983 289-299

\art RV J. Raissy, L. Vivas: Dynamics of two-resonant biholomorphisms! Math. Res. Lett.! 20 2013 757-771

\art R J. Rib\'on: Families of diffeomorphisms without periodic curves! Michigan Math. J.! 53 2005 243-256

\book Ri1 M. Rivi: Local behaviour of discrete dynamical systems! Ph.D. Thesis,
Universit\`a di Firenze, 1999

\art Ri2 M. Rivi: Parabolic manifolds for semi-attractive holomorphic germs! Michigan Math. J.! 49 2001 211-241

\art Ro1 F. Rong: Quasi-parabolic analytic transformations of $\C^n$! J. Math. Anal. Appl.! 343 2008 99-109

\art Ro2 F. Rong: Quasi-parabolic analytic transformations of $\C^n$. Parabolic manifolds! Ark. Mat.! 48 2010 361-370

\art Ro3 F. Rong: Robust parabolic curves in $\C^m$ ($m\ge3$)! Houston J. Math.! 36 2010 147-155

\art Ro4 F. Rong: Absolutely isolated singularities of holomorphic maps of $\C^n$ tangent to the identity! Pacific J. Math.! 246 2010 421-433

\art Ro5 F. Rong: Parabolic manifolds for semi-attractive analytic transformations of $\C^n$! Trans. Amer. Math. Soc.! 363 2011 5207-5222

\art Ro6 F. Rong: Local dynamics of holomorphic maps in $\C^2$ with a Jordan fixed point! Michigan Math. J.! 62 2013 843-856

\art Ro7 F. Rong: Reduction of singularities of holomorphic maps of $\C^2$ tangent to the identity! Publ. Math. Debrecen! 83/4 2013 537-546

\art Ro8 F. Rong: The non-dicritical order and attracting domains of holomorphic maps tangent to the identity! Int. J. Math.! {25 \rm 1450003} 2014 1-10

\pre Ro9 F. Rong: New invariants and attracting domains for holomorphic maps in $\C^2$ tangent to the identity! Preprint, to appear in Publ. Mat.! 2014

\pre S D. Sauzin: Introduction to $1$-summability and resurgence! Preprint, arXiv:1405.0356! 2014

\art Sh A.A. Shcherbakov: Topological classification of germs of conformal
mappings with identity linear part! Moscow Univ. Math. Bull.! 37 1982 60-65

\art SV B. Stensones, L. Vivas: Basins of attraction of automorphisms in $\C^3$! Ergodic Theory Dynam. Systems! 34 2014 689-692

\art U1 T. Ueda: Local structure of analytic transformations of two complex variables, I! J. Math. Kyoto Univ.! 26 1986 233-261

\art U2 T. Ueda: Local structure of analytic transformations of two complex variables, I\negthinspace I! J. Math. Kyoto Univ.! 31 1991 695-711

\art Us S. Ushiki: Parabolic fixed points of two-dimensional complex dynamical
systems! S\=uri\-kaise\-kiken\-ky\=usho K\=oky\=uroku! 959 1996 168-180

\art V1 L. Vivas: Degenerate characteristic directions for maps tangent to the identity! Indiana Univ. Math. J.! 61 2012 2019-2040

\art V2 L. R. Vivas: Fatou-Bieberbach domains as basins of attraction of automorphisms 
tangent to the identity! J. Geom. Anal.! 22 2012 352-382 

\art W B. Weickert: Attracting basins for automorphisms of $\C^2$! Invent. Math.! 132 1998 581-605

\bye